\documentclass[11pt]{article}
\usepackage{color, amsmath, amssymb, amsbsy, amsthm, graphicx, bbm, amsfonts, bm, dsfont}
\usepackage{setspace, booktabs, comment, graphicx, algorithm, algorithmic, latexsym}
\usepackage[top=1in, bottom=1in, left=1in, right=1in]{geometry}
\usepackage[toc,page]{appendix}
\usepackage[colorlinks,allcolors=blue]{hyperref}
\allowdisplaybreaks

\RequirePackage[authoryear]{natbib}
\theoremstyle{plain}
\newtheorem{theorem}{Theorem}[section]

\newtheorem{lemma}{Lemma}[section]

\newtheorem{proposition}{Proposition}[section]

\def \hat{\widehat}
\def \tilde{\widetilde}

\usepackage{tikz}
\usetikzlibrary{shapes,decorations,arrows,calc,arrows.meta,fit,positioning}
\tikzset{
	-Latex,auto,node distance =1 cm and 1 cm,semithick,
	state/.style ={ellipse, draw, minimum width = 0.7 cm},
	point/.style = {circle, draw, inner sep=0.04cm,fill,node contents={}},
	bidirected/.style={Latex-Latex,dashed},
	el/.style = {inner sep=2pt, align=left, sloped}
}

\usepackage[all,arc]{xy}
\usepackage{enumerate}
\usepackage{mathrsfs}

\usepackage{amsmath}
\usepackage{amsthm}
\usepackage{amsfonts}
\usepackage{amssymb}
\usepackage{amsbsy}
\usepackage{graphicx}
\usepackage{enumitem}
\usepackage{physics}
\usepackage{hyperref}
\usepackage{bbm}

\newcommand{\mc}[1]{\mathcal{#1}}

\newcommand{\R}{\mathbb{R}}

\newcommand{\E}{\mathbb{E}}

\newcommand{\ra}{\rightarrow}
\newcommand{\toinf}{\ra \infty}

\newcommand{\beq}{\begin{equation}}
\newcommand{\eeq}{\end{equation}}

\theoremstyle{definition}

\newcommand{\Var}[1]{\mathrm{Var}(#1)}
\newcommand{\lra}{\longrightarrow}
\numberwithin{equation}{section}
\numberwithin{theorem}{section}
\newcommand{\ind}{\mathbbm{1}}
\newcommand{\ptl}{\partial}

\DeclareMathOperator*{\argmin}{arg\,min}

\newcommand{\chacorr}{C}
\newcommand{\renyicorr}{R}

\newcommand{\renyi}{R\'{e}nyi}

\newcommand{\dens}{\phi}

\title{Measures of independence and functional dependence}
\author{Peter J. Bickel \\ University of California, Berkeley}

\begin{document}


\maketitle

\begin{abstract}
We follow up on Shi et al's (2020) and Cao's and my  (2020) work on the local power of a new test for independence, Chatterjee (2019), and its relation to the local power properties of classical tests. We show quite generally that for testing independence with local alternatives either Chatterjee's rank test has no power or it may be misleading: The Blum,Kiefer,Rosenblatt and other omnibus classical rank tests do have some local power in any direction other than those where significant results may be misleading. We also suggest methods of selective inference in independence testing. Chatterjee's statistics like Renyi's \cite{RenyiMD1959} also identified functional dependence. We exhibit statistics which have better power properties than Chatterjee's but also identify functional dependence.
\end{abstract}

\section{Introduction}\label{section:introduction}

In a remarkable paper \cite{ChatterjeeNCCA2020}, Sourav Chatterjee proposed a new coefficient of correlation based on an i.i.d. sample $(X_i, Y_i), i = 1, \ldots, n$. Assuming there are no ties among the $X_i$'s and $Y_i$'s (see \cite{ChatterjeeNCCA2020} for the definition in the general case), the correlation is defined as
\beq\label{eq:chacorr-est-def} \hat{\chacorr}_n(X, Y) := 1 - \frac{3 \sum_{i=1}^{n-1} |r_{i+1} - r_i|}{n^2-1}, \eeq
where the $r_i$ are defined as follows. First, sort $X_{(1)} \leq \cdots \leq X_{(n)}$, and for each $i$ let $Y_{(i)}$ be the $Y$ sample corresponding to $X_{(i)}$. Then $r_i$ is defined as the rank of $Y_{(i)}$, i.e. the number of $j$ such that $Y_j \leq Y_{(i)}$. Chatterjee showed that as $n \toinf$, $\hat{\chacorr}_n$ a.s. converges to the population measure
\begin{align*}
    \label{eq:chacorr-def} \chacorr(X, Y) := \frac{\int \Var{\E[\ind(Y \geq t) ~|~ X]} d\mu(t)}{\int \Var{\ind(Y \geq t)} d\mu(t)}
\end{align*}

\begin{align}
 = 1 - \frac{\int \E [\mathrm{Var}(\ind(Y \geq t) ~|~ X) ] d\mu(t)}{\int \Var{\ind(Y \geq t)} d\mu(t)},
\end{align}

where $\mu$ is the law of $Y$. Here $Y$ is assumed to not be constant. In the case where $X, Y$ are continuously distributed, this measure was introduced by Dette et al. \cite{DETTECNMR2013}. The measure $\chacorr$ has a number of interesting properties:
\begin{enumerate}[label=\Alph*)]
    \item $0 \leq \chacorr \leq 1$.
    \item $\chacorr = 0$ if and only if $X$ and $Y$ are independent.
    \item $\chacorr = 1$ if and only if $Y = h(X)$ a.s. for some measurable function $h : \R \ra \R$.
    \item $\chacorr$ is asymmetric, but can be easily symmetrized to
    \[ \chacorr^*(X, Y) := \max(\chacorr(X, Y), \chacorr(Y, X)),\]
    which clearly satisfies $\chacorr^* = 1$ if and only if $X$ is a function of $Y$ or $Y$ is a function of $X$ (or both).
    \item $\chacorr$ is invariant under strictly increasing transformations of $X$ and $Y$ separately.
\end{enumerate}

This measure is akin to the {\renyi} correlation (also commonly called the maximal correlation), which we shall denote $\renyicorr$ or $\renyicorr(X, Y)$, and is defined as the maximum Pearson correlation between all pairs of $L^2$ functions of $X$ and $Y$ respectively. $R$ may be computed as the square root of the maximal eigenvalue of a compact self adjoint operator,
\[ T : L^2_0(X) \ra L^2_0(X), \]
(or $L^2_0(Y) \ra L^2_0(Y)$ with appropriate changes), where $L^2_0(X)$ is the subspace of $L^2(X)$ consisting of mean zero random variables. The operator $T$ is given by
\beq T(f(X)) := \E[\E[f(X) ~|~ Y] ~|~ X]. \eeq
The {\renyi} correlation is well known to have properties $A, B,$ and $C$, but it is symmetric, and 
\begin{enumerate}[label=C*)]
    \item $\renyicorr = 1$ if and only if $g(Y) = h(X)$ for some functions $g$ and $h$, with $g(X) \in L^2(X)$ and $h(Y) \in L^2(Y)$.
\end{enumerate}
An extensive account of the history and  computation of estimates of R, and other properties of $\renyicorr$ may be found for instance in \cite{BickelEAES1998, RenyiMD1959}.

It is clear that estimators of $R$, generically call them $\hat{R}$, involve implicitly estimation of the conditional distribution of $Y$ given $X$, and conversely. Thus, except in the finite discrete case, they require assumptions on smoothness of the joint distribution  of X and Y  for convergence to $R$. In testing the hypothesis $H: \textit{X independent of Y}$ versus contiguous alternatives, $\hat{R}$ has no power, i.e. has asymptotic power equal to the level of significance. 

Perhaps surprisingly, it has been shown independently by Shi et al.(2020) \cite{ShiPCRCA2020}  and Cao and Bickel \cite{CaoCTEPO2020} (2020) that the test based on $\hat{\chacorr}_n$ and the related one based on Dette et al's statistic have the same poor behaviour as tests of independence in many copula models. 

In this paper, we show quite generally that Chatterjee's test has no power (power $\alpha$) against local parametric alternatives (those converging to $H$ at rate $O(\frac{1}{\sqrt{n}})$ ) other than a special set of parametric families. Those alternatives can yield significant results but are equivalent in a strong sense to models where $X$ and $Y$ are always independent, varying only in their marginal distributions. 

We go on to show that this behaviour is to be expected from a large family of rank tests which use the ranks of $X$'s belonging to adjacent $X$'s. Or, put another way, these tests look at departures from independence from a very local point of views in $(X,Y)$ space.

Following an approach of Bickel, Ritov,Stoker (2006) \cite{BickelTTGFA2006}, we also introduce a class of simple tests, computable in $O(n)$ steps, which have power tailored to the complexity of the alternative, and develop this into an approach to selective inference.

In section 3, we discuss functional dependence and a principle underlying both Renyi's and Chatterjee's statistics, deriving some new procedures and relating these to classical ones in a general context, including multivariate and discrete ones.

Finally, in Section 4, we discuss the various ways by which parameters satisfying both of Chatterjee's properties B and C can be constructed with satisfactory behaviour under both regimes.

\section{Tests of independence}
\subsection{Definitions}

There has recently been a surge of interest in such tests. From the point of view of minimax analysis (\cite{KimMOPTM2020}), local power analysis (\cite{ShiPCRCA2020,WeihsSRCGS2018}, \cite{CaoCTEPO2020}, many types of statistics (\cite{GenestTIADM2019,LinBPCRA2021,SzekelyMTDCD2007}) have been proposed.

In addition to Chatterjee’s and Dette et al’s  $C$, there are a number of well known population parameters corresponding to test statistics based on ranks or more generally  invariant under monotone transformations of coordinates which are consistent against essentially all fixed alternatives to $H:$ $X$ independent of $Y$. Several prominent ones are due to Hoeffding(1948), Blum,Kiefer,Rosenblatt(1961) \cite{BlumDFTIJ1961}, Bergsma and Dassios(2014) \cite{BergsmaCTIBM2014}, and more recently Weihs, Drton, Meinshausen(2018) \cite{WeihsSRCGS2018} and many others (Kim et al(2020) \cite{KimMOPTM2020}). In addition there are even more classical parameters corresponding to directed statistics such as Spearman’s $\rho$ and Kendall’s $\tau$. All of these statistics, $T_n$, can be related to population parameters $T(F)$ where $F$ ranges over all bivariate distributions, represented as cdfs, and which have the property that $T(F)=0$ if $X$ and $Y$ are independent. This is Weihs et als \cite{WeihsSRCGS2018}, property I. The statistics $T_n$ are typically, but not necessarily, obtained by plug in but will be identified by consistency to the population parameters as $n \to \infty$. Plug in means,
\begin{align*}
T_n \equiv T(\hat{F}) 
\end{align*}
,where $\hat{F}$  is the empirical cdf  of the sample $(X_i,Y_i)$ from $F$.

The class of parameters we consider are functions, not just of $F$,
but of $C$, defined as the cdf of the vector $(U,V)$ with  $U=F_1(X)$, $V=F_2(Y)$, where $F_1,F_2$ are the cdf of $X,Y$ respectively, 
and $U$ and $V$ have corresponding marginals $C_1$ and $C_2$. If $F_1$ and $F_2$ are continuous, C is a copula and $C_1, C_2$ are both $U(0,1)$. In general,
\begin{align*}
C(u,v)=F(Q_1(u),Q_2(v))\\
\end{align*}
,where $Q_j(.)$ is the inverse of $F_j$ and the inverse is defined for any cdf $F$ by,     
\begin{align*}
Q(u)=\inf\{x:F(x) \geq u\}\\
\end{align*}

The corresponding $T_n$ will depend on the empirical cdf of the  empirical probability integral transforms, $\hat{F_1}(X),\hat{F_2}(Y))$. This object, we shall call $\hat{C}$, is referred to as the empirical copula \cite{TsukaharaSECM2005}. Thus,
\begin{align}
\hat{C}(u,v)=\hat{F}(\hat{Q_1}(u) ,\hat{Q_2}(v))
\end{align}
 where $\hat{Q_j}$ are just the quantile functions of $\hat{F_j}$.

Up to this point we have considered all $F$, continuous, discrete or both. For testing $H$, the assumption that $F_j,j=1,2$ are continuous, manifested by no ties in the $X$'s or $Y$'s, simplifies matters considerably, as we have partly noted. The distribution of $T(\hat{C})$ under $H$ does not depend on F. It is always a function of $C$, the cdf of the $(U_i, V_i), i=1,...,n$ only.

If $\hat{C}$ is the empirical cdf of the vectors of ranks normalized by n, being a function of $\hat{C}$ simply means that the $T_n$ are ordinary rank statistics. That is,
$T_n=T_n(S_{(1)},…,S_{(n)})$ where  $S_{(j)}$ is the rank of the $Y$ corresponding to the order statistic $X_{(j)}$. All the classical statistics for testing independence, as well as the more recent ones of Weihs et al., Bergsma et al., as well as Chatterjee and Dette et al., are functions of the empirical copula.

As examples, we consider,

1) The Blum Kiefer-Rosenblatt parameter

\begin{align}\label{T_BKR}
T_{BKR}(F) &\equiv \int_{-\infty}^{\infty}\int_{-\infty}^{\infty}(F(x,y)-F_{1}(x)F_{2}(y))^2dF_{1}(x)dF_{2}(y) \notag\\
&= \int_{0}^{1}\int_{0}^{1}(C(u,v)-C_1(u)C_2(v))^2dudv \notag\\
&= \int_{0}^{1}\int_{0}^{1}(C(u,v)-uv)^2dudv
\end{align}
in the copula case.

2) The Spearman correlation parameter

\begin{align}
    T_{S}(F)&\equiv corr(F_{1}(x),F_{2}(Y))\\
         &=\frac{\int_{0}^{1}\int_{0}^{1}(uv-\mu^{*}_{U}\mu^{*}_{V})dC(u,v)}{\sigma^{*}_{U}\sigma^{*}_{V}} \notag
\end{align}
where $\mu^{*}_{U}, \mu^{*}_{V}, \sigma^{*}_{U}, \sigma^{*}_{V}$ are the means and variances of $U\equiv F_{1}(X), V \equiv F_{2}(Y)$.

If C is a copula, $T_{S}(F)=12[\int_{0}^{1}\int_{0}^{1} uv dC(u,v)-\frac{1}{4}]$

3) The sup norm parameter, we call $T_K$ for Kolmogorov, introduced in \cite{BlumDFTIJ1961}.

\begin{align*}
    T_{K}(F) &\equiv \sup\{|F(x,y)-F_1(x)F_2(y)|: x,y\} \\
    &= \sup \{C(u,v)-C_1(u)C_2(v): u,v \in [0,1]^2\}
\end{align*}

The familiar corresponding statistics are obtained up to $O_{p}(\frac{1}{n})$ by replacing $F$ by $\hat{F}$ or $C$ by $\hat{C}$. If $C$ is a copula, parameters and statistics simplify since $C_1(v)=C_2(v)=v$, and statistics are functions of the normalized ranks.

The statistics can then be written in terms of $r_i$ defined for (1.1) as,

\begin{align*}
    \hat{T}_{S} &= \frac{12}{n} \sum_{i=1}^{n}(\frac{ir_i}{n^2}-\frac{1}{4})+O_p(\frac{1}{n})\\
    \hat{T}_{K} &= max\{|\frac{1}{n}\sum_{i=1}^{j}(1(r_i \leq k)-\frac{jk}{n^2})| : 1\leq j,k \leq n\} +O_p(\frac{1}{n})\\
     \hat{T}_{BKR} &= \frac{1}{n^2} \sum_{i,j=1}^{n}A_{ij} B_{ij}-
     \frac{1}{2n} \sum_{i=1}^{n}(\frac{ir_i}{n})^2-\frac{7}{18}+O_p(\frac{1}{n})
\end{align*}
where $A_{ij}=\frac{i}{n} \wedge \frac{j}{n}$, $B_{ij}=\frac{r_i}{n} \wedge \frac{r_j}{n}$.

Both $\hat{T}_{BKR}$ and $\hat{T}_{K}$ require $O(n^2)$ operations. $\hat{T}_{S}$ requires $O(n\log n)$ operations.

There are, of course, many other parameters measuring independence, discussed in Tsukahara (2005), Weihs et al (2018) among others. The corresponding statistics are typically plug in. Most, such as $\hat{T}_{S}, \hat{C}_n$, are asymptotically Gaussian, and linear, in the sense that,

\begin{align}
    T(\hat{F})&= T(F) + \frac{1}{n} \sum_{i=1}^{n}{\Psi(X_i,Y_i,F)}+o_{p}(n^{-\frac{1}{2}})
\end{align}

with influence function $\Psi \in L_{2}(F), E_{F} \Psi(X,Y,F)=0$.\\
    
Some are not, for example $\hat{T}_{BKR}$ and $\hat{T}_{K}$, given by, 

\begin{align*}
    \hat{T}_K \equiv sup \{ |\hat{C}(u,v)-\hat{C_1}(u)\hat{C_2}(v)|:u,v\}.
\end{align*}

The limiting behavior of these statistics is formally derivable using the stochastic process\\ $\sqrt{n}(\hat{C}-C)$ which we are about to introduce. In some cases, such as the Weihs et al statistics. $U$ statistic theory provides much simpler proofs. By (1.2), the Chatterjee-Dette et al parameter, which we shall call $T_{C}(F)$ ,can also be expressed in terms of $F$,
\begin{align}
    T_{C}(F)&=\frac{ E_{F}\bigg(Var_{F}(1(Y\leq Y')|X)\bigg)}{Var_{F}(1(Y\leq Y'))}
\end{align}

where $(X,Y) \sim F$ and $Y'$ is independent of $(X,Y)$ and distributed according to $F_2(.)$.\\ 

We return to this definition in section 3.4. $C$ can replace $F$ in this definition. However, the corresponding statistics in \cite{ChatterjeeNCCA2020} and \cite{DETTECNMR2013} are not plug ins.  Our formulae so far are all valid for general $F$, $\hat{F}$. For the rest of this section, we focus on the copula case, F continuous. We return to the general case at the end of the section, in a brief discussion.

\subsection{Stochastic Processes associated with the hypothesis $H$ of independence}

The process,

\begin{align}
    T_n(u,v)=\sqrt{n}\bigg(\hat{C}(u,v)-C(u,v)\bigg)
\end{align}

has been heavily studied, (\cite{GenestTIADM2019},\cite{TsukaharaSECM2005}, for example). In general, its behavior requires delicate conditions because of its dependence on $\hat{Q}_{j}$. However, for the situation we focus on, $X,Y$ independent and $F$ continuous, the following special case of the general theory characterizing $T_n(.,.)$ is simple.

Let $E_n, E_{n1}$ and $E_{n2}$ be the empirical process of $\{(U_i,V_i)\}, \{U_i\}$ and $\{V_i\}$ respectively where $U_i=F_1(X_i), V_i=F_2(Y_i)$ are i.i,d. $U(0,1)$. Then, $C(u,v)=uv$ and $\hat{C}(.,.)$ is the same as the empirical copula of the $(U_i,V_i)$. Thus, if we simplify notation and set $F_1,F_2$ to the identity,

\begin{align}
    \sqrt{n}\bigg(\hat{C}(u,v)-C(u,v)\bigg)=E_n(\hat{Q}_1(u),\hat{Q}_2(v)
    )+\sqrt{n}(\hat{Q}_1(u)\hat{Q}_2(v)-uv).
\end{align}

By classical empirical process theory and Bahadur \cite{BahadurNQLS1966},

\begin{align}
    \sqrt{n}(\hat{Q}_j(u)-u)=-E_{nj}(u)+o_p(1), j=1,2
\end{align}

and,

\begin{align}
    E_n(Q_1(u),Q_2(v))=E_n(u,v)+o_p(1).
\end{align}

Note that $o_p(1)$ and similar expressions denote the order of functions in the $L_{\infty}$ norm.

By the delta method,

\begin{align}
    \sqrt{n}(\hat{Q}_1(u)\hat{Q}_2(v)-uv)=u\sqrt{n}(\hat{Q}_2(v)-v)+v\sqrt{n}(\hat{Q}_1(u)-u)\\
    =-uE_{n2}(v)-vE_{n1}(u)+O_p(n^{-\frac{1}{2}})
\end{align}

by (2.8).\\

We conclude that,

\begin{proposition}
If $X$ and $Y$ are independent and $F$ is continuous,
$T_n(u,v)=E_n(u,v)-uE_{n2}(v)-vE_{n1}(u)+o_p(1)$
\end{proposition}

This is a special case of Proposition 1 of Tsukahara (2005)\cite{TsukaharaSECM2005}.

We now claim that, under the same conditions, 

\begin{align}
    \tilde{\Lambda}_n(u,v) \equiv uE_{n2}(v)+vE_{n1}(u)=\Pi(E_n(u,v)|\Lambda_n)
\end{align}

Here $\Pi(.|\Lambda_n)$ denotes projection in $L_2(U_i,V_i: i=1,...,n)$ on the linear subspace $\Lambda_n$ given by,

\begin{align}
    \Lambda_n=\{\sum_{i=1}^{n}{(f_i(U_i)+g_i(V_i)}: Ef_{i}^2(U) < \infty,  Eg_{j}^2(V) < \infty, Ef_{i}(U)= Eg_{j}(V)=0 \textit{ for all } i,j\}
\end{align}

This follows from the well known observations:

\begin{enumerate}
    \item If $h(U,V) \in L_2(U,V)$, the projection of $h$ on the space of all functions of $U$ is $E(h(U,V)|U)$ and similarly for functions of $V$
    \item The projection operator is additive on a sum of orthogonal subspaces.
\end{enumerate}

Thus, our claim is equivalent to,

\begin{align*}
    E(1(U\leq u,V \leq v)|V)-uv=u(1(V\leq v)-v)
\end{align*}

and similarly for $U$, and (2.12) follows.

Denote the process of (2.12) by $\Lambda_n(u,v)$. Let $D[0,1]^2$ be defined as in Billingsley \cite{BillingsleyCPM1968}.

\begin{proposition}
Under the conditions of Proposition 2.1, let $g,h:D[0,1]^2 \to D[0,1]^2, S_n \equiv g(T_n(.,.)), V_n \equiv h(\tilde{\Lambda}_n(.,.))$. $g,h$ measurable on $D[0,1]^2$.\\
Then, if $(S_n,V_n) \xrightarrow{d} (S,V)$, $S$ and $V$ are independent.
\end{proposition}

$\bold{Proof:}$ From classical empirical process theory, $(T_n(.,.),\tilde{\Lambda}_n(.,.))$ converge jointly weakly to jointly Gaussian processes with the same covariance structure as $E_n(u,v)-uE_{n2}(v)-vE_{n1}(u)$ and $uE_{n2}(v)+vE_{n1}(u)$, which we shall call $T^*$ and $\Lambda$. We shall denote the limit of $E_n(.)$ by $E(.)$ and those of $E_{nj}(.)$ by $E_j(.)$. Note that we can think of these limit processes as a.s. limits in the $l_\infty$ norm. By (2.12), $\tilde{\Lambda}_n$ and $E_n-\tilde{\Lambda}_n$ are orthogonal. Their Gaussian limits must then be independent.\qedsymbol

\newcommand{\tend}[1]{\hbox{\oalign{$\bm{#1}$\crcr\hidewidth$\scriptscriptstyle\bm{\sim}$\hidewidth}}}

Formulae for the covariances functions of $\tilde{\Lambda}_n$ and $T_n$ and hence of $T^*$ follow. If $\tend{u}\equiv (u_1,u_2)$, $\tend{v}\equiv (v_1,v_2)$,

\begin{align}
    cov(T^*(\tend{u}),T^*(\tend{v})) &= cov(E_n(\tend{u}),E_n(\tend{v}))-cov(\tilde{\Lambda}_n(\tend{u}),\tilde{\Lambda}_n(\tend{v})) \notag\\
    &= (u_1 \wedge u_2)(v_1 \wedge v_2)-u_1u_2v_1v_2-u_1v_1[(u_2\wedge v_2)-u_2v_2]-u_2v_2[(u_1\wedge v_1 -u_1v_1)] \notag\\
    &=(u_1 \wedge v_1-u_1v_1)(u_2 \wedge v_2-u_2v_2)
\end{align}
as noted by Tsukahara. We can represent

\begin{align*}
\Lambda(u,v)=uW^{\circ}(v) +vW^{\circ}(u)    
\end{align*}

where $W^{\circ}$ is a Brownian Bridge.

We proceed with the implication of this result for local power.

\subsection{Local power against H}

\subsubsection{Le Cam Theory}

We make the following assumptions,following Le Cam and Yang (1990) \cite{CamASBC1990},

\begin{enumerate}[label=\Roman*]
    \item We begin with a model $\mc{M}=\{h_\theta(x, y): |\theta| < 1 \}$, where $h_\theta(x, y)$ is a joint density with respect to Lebesgue measure, and $h_0(x, y) = f_0(x) g_0(y)$ (independence). 
    \item Suppose, the score function,
    \[ \dot{\ell}_0(x, y) := \frac{\partial}{\partial \theta} \log h_\theta(x, y) \bigg|_{\theta = 0} \]
    exists, and suppose the family $\{h_\theta(x, y)\}_{|\theta| < 1}$ is quadratic mean differentiable (QMD) at $\theta = 0$ with score function $\dot{\ell}_0$. That is, we have
    \[ \int \int \big(\sqrt{h_\theta(x, y)} - \sqrt{h_0(x, y)} - (1/2) \theta \dot{\ell}_0(x, y) \sqrt{h_0(x, y)}\big)^2 dx dy = o(\theta^2). \]
    Moreover, assume $\E_0 \dot{\ell}_0(X, Y)^2 > 0$ (the assumption of quadratic mean differentiability implies $\E_0[ \dot{\ell}_0(X, Y)^2 ]< \infty$ and $\E_0 \dot{\ell}_0(X, Y) = 0$). 
\end{enumerate}

For $\theta \in \R$ let $P_{n, \theta}$ be the product measure on $([0, 1] \times [0, 1])^n$ with density $\prod_{i=1}^n h_\theta(x, y)$. Take $t \in \R$ and let $\theta_n := t / \sqrt{n}$. By Le Cam's Theorem (see e.g. \cite[Theorem 7.2 and Example 6.5]{VaartAS1998}), under our assumptions $\{P_{n, \theta_n}\}_{n \geq 1}$ is contiguous with respect to $\{P_{n, 0}\}_{n \geq 1}$. That is if $Q_n$ is a function of $((X_i, Y_i), 1 \leq i \leq n)$, then
\[ Q_n \stackrel{P_{n, 0}}{\lra} 0 \text{ implies }  Q_n \stackrel{P_{n, \theta_n}}{\lra} 0. \]
For our purposes, more importantly, Le Cam's third lemma holds (see e.g. \cite[Theorem 7.2 and Example 6.7]{VaartAS1998}), stating that if, 

\begin{align}\label{eq:le-cam-3}
    \bigg(\sqrt{n}Q_n, \frac{1}{\sqrt{n}} \sum_{i=1}^n \dot{\ell}_0(X_i, Y_i)\bigg) \stackrel{d}{\lra} N\bigg(\begin{pmatrix} 0 \\ 0 \end{pmatrix}, \begin{pmatrix} \sigma^2 & \rho \sigma \tau_0 \\ \rho \sigma \tau_0 & \tau_0^2 \end{pmatrix}\bigg) 
\end{align}

under $\{P_{n, 0}\}_{n \geq 1}$, where $\tau_0^2 := \E_0 [\dot{\ell}_0(X, Y)^2]$ is the Fisher information, then under $\{P_{n, \theta_n}\}_{n \geq 1}$, we have
\[ \sqrt{n} Q_n \stackrel{d}{\lra} N(tc, \sigma^2), \]
where $c := \rho \sigma \tau_0$ is the asymptotic covariance (under $\{P_{n, 0}\}_{n \geq 1}$) between the statistic $T_n$ and the score statistic. 

Assuming \eqref{eq:le-cam-3} holds under $\{P_{n, 0}\}_{n \geq 1}$, $\sqrt{n} T_n / \sigma$ can be viewed as a test statistic for the hypothesis of independence, while

\begin{align}
 L_n := \frac{1}{\tau_0 \sqrt{n}} \sum_{i=1}^n \dot{\ell}_0(X_i, Y_i)  
\end{align}

can be viewed as the asymptotically optimal test statistic for the family $\{h_\theta(x, y)\}_{|\theta| < 1}$. The Pitman efficiency (see e.g. \cite[Chapter 8]{VaartAS1998}) of the first test to the second is by (2.16), 
\[ e(Q_n) = \rho^2. \]
Note that $e(L_n) = 1$.

$\bold{Remark:}$ 1) Note that if $\dot{\ell}_0(x, y)$ exists as a pointwise derivative and is continuous in $L_2$ norm at $\theta=0$,  QMD follows.

2) Let $U \equiv F_1(X), V \equiv F_2(Y)$ where $F_1, F_2$ are the cdf of $X$ and $Y$ respectively if $\theta=0$. If we equivalently view $\mathcal{M}$ as the model of distributions of $(U,V)$, it is easy to see that the corresponding score function is unchanged as a random variable and is given by $\dot{l}_0(F_1^{-1}(U),F_2^{-1}(V))$. In what follows, we shall usually WLOG assume $F_1,F_2$ are the identity and hence $X=U$, $Y=V$ both $\mathcal{U}(0,1)$ under $H$. We equivalently use $\dot{l}_0(U,V)$ for $\dot{l}_0(X,Y)$.

We will use a slight generalization of the third lemma, which is evident once we note that contiguity preserves tightness.

\begin{lemma}
Let $M_n(.)$ be a stochastic process taking values in $D^{k} \equiv D([0,1]^{k})$ (D as defined in Billingsley \cite{BillingsleyCPM1968} ). Suppose under $\{P_{n,0}\}_{n\geq 1}$,
    \begin{align*}
        \bigg(\sqrt{n}M_n(.), \frac{1}{\sqrt{n}} \sum_{i=1}^n \dot{\ell}_0(X_i, Y_i)\bigg) \stackrel{d}{\lra}  
        \bigg( M(.), Z \bigg)
    \end{align*}
where $\stackrel{d}{\lra}$ indicates weak convergence and $\bigg( M(.), Z \bigg)$ is a Gaussian process, $\mathbb{E}M(.)=0$, and $cov \bigg( M(s), M(t) \bigg) = \sigma (s,t) \textit{ for } \sigma: D^{k}\times D^{k} {\lra} \R$.

Let,
\begin{align*}
    C(s)=cov \bigg( M(s), Z \bigg)
\end{align*}
Then, under $\{P_{n,\theta_{n}}\}$,
\begin{align*}
    \sqrt{n}M_n(.) \stackrel{d}{\lra}  M(.)+tC(.)
\end{align*}
\end{lemma}

\subsubsection{Some Definitions}

Suppose $T_n$ is a rank statistic. and a critical region is $[T_n \geq t_n]$. A test based on $T_n$ has \emph{asymptotic level} $\alpha$ iff
\begin{align*}
    \lim_{n} P_0[T_n \geq t_n]=\alpha.
\end{align*}

A level $\alpha$ test is \emph{asymptotically unbiased} iff for all regular one dimensional parametric submodels with efficient score function $l_{0}^{*}$ and $|\theta_n| \leq \frac{M}{\sqrt{n}}, M<\infty$,

\begin{align*}
    \liminf P_{\theta_n}[T_n \geq t_n]>\alpha.
\end{align*}

It has \emph{"no power"} asymptotically if for all submodels as above, $\frac{\epsilon}{\sqrt{n}} \leq |\theta_n| \leq \frac{M}{\sqrt{n}},$ all $\epsilon>0, M<\infty$,

\begin{align}
    \liminf P_{\theta_n}[T_n \geq t_n]=\alpha.
\end{align}

\subsubsection{Local power behaviour of $\hat{C}_n$}
We apply our machinery to $\hat{C}_n$.

We note, if $H$ is true and $F$ is continuous, then by Angus(1995), see also Shi et al (2021) and Cao and Bickel (2021), we have
\begin{align}
    \hat{C}_n(X,Y)=\frac{1}{n} \sum_{i=1}^{n}{(\frac{2}{3}-|V_{[i]}-V_{[i+1]}|-2V_{[i]}(1-V_{[i]}))}+o_p(n^{-\frac{1}{2}})
\end{align} 

where we define $V_{[i]}$ by, $V_{[i]}=V_j$ iff $Rank(U_j)=i$, and $U_i \equiv F_1(X_i),V_i \equiv F_2(Y_i)$.\\

We claim that,

\begin{proposition}
If $\dot{l}_0$ is as defined in (2.16) and 
\begin{align}
    E_0(\dot{l}_0(X,Y)|X)=E_0(\dot{l}_0(X,Y)|Y)=0
\end{align}
or equivalently, 
\begin{align}
    \dot{l}_0(U_1,V_1) \in \Lambda_1^{\perp}
\end{align}
then $\hat{C}_n$ has no power.
\end{proposition}

\textbf{Proof: } This proposition follows from the more general Theorem 2.3 and the following.

Let $\mathcal{M}$ be a model as defined in Section 2.3.1 with density $h_{\theta}(.,.)$ and score function $\dot{l}_0(X,Y)$. Without loss of generality we let $F_1=F_2=Identity$ and write $\dot{l}_0(U,V)$. Decompose $\dot{l}_0$ as, 
\begin{align*}
    \dot{l}_0(U,V)=\dot{l}_{00}(U,V)+\dot{l}_{01}(U,V)
\end{align*}
where $\dot{l}_{00}\in \Lambda_1$ and $\dot{l}_{01}\in \Lambda_1^{\perp}$,
where $\Lambda_1$ is defined in (2.13).

For a statistic $Q_n$, satisfying the hypothesis of Section 2.3.1, let,

\begin{align*}
    \lambda_0(Q_n) &\equiv lim_n cov_0(\frac{1}{\sqrt{n}}\sum_{i=1}^n \dot{l}_{00}(U_i,V_i),Q_n)\\
    \lambda_1(Q_n) &\equiv lim_n cov_0(\frac{1}{\sqrt{n}}\sum_{i=1}^n \dot{l}_{01}(U_i,V_i),Q_n)
\end{align*}

If $L_n$ as defined by (2.16) and $Q_n$ w.l.o.g. both have asymptotic variance $1$, the Pitman efficiency of $Q_n$ is given by,
\begin{align*}
    e(Q_n)=(\lambda_0(Q_n)+\lambda_1(Q_n))^2
\end{align*}

We shall show in Theorem 2.3 that $\lambda_1(\hat{C}_n)=0$.

Proposition 2.3 follows since its assumption is that $\dot{l}_{00}=0$ and hence $\lambda_0(\hat{C}_n)=0$.

Going further, write,
\begin{align*}
    \dot{l}_{00}(U,V)=a(U)+b(V)
\end{align*}
with $E_0 a(U)=E_0 b(V)=0$.

Suppose,
\begin{align*}
    h_0(x,y)=f_0(x)g_0(y).
\end{align*}

Define a new model $\tilde{\mathcal{M}}$ by,

\begin{align*}
    \tilde{\mathcal{M}}=\{\tilde{h}_{\theta}(x,y): |\theta| <1 \}
\end{align*}

and
\begin{align*}
    \tilde{h}_0(x,y)=h_0(x,y).
\end{align*}

Define,
\begin{align*}
    \tilde{h}_{\theta}(x,y)= \frac{h_0(x,y)e^{\theta(a(x)+b(y))}}{Z(\theta)}.
\end{align*}
where $Z(0)=1$ and $Z(\theta)$ is a normalizing function.

Then, $X$ and $Y$, or equivalently, $U$ and $V$ are independent for all $\theta$.

Let,

\begin{align}
    \Delta_n(\theta) \equiv \int...\int |\Pi_{i=1}^n h_{\theta}(x_i,y_i)-\Pi_{i=1}^n \tilde{h}_{\theta}(x_i,y_i)| \Pi_{i=1}^ndx_idy_i
\end{align}

We will, show under some regularity conditions,

\begin{align}
    sup\{\Delta_n(\theta): |\theta| \leq \frac{M}{\sqrt{n}}\} \to 0
\end{align}

We prove the following in the appendix.

\begin{proposition}
Let $h_{\theta}, \tilde{h}_{\theta}, a, b$ be as given.

Assume, in addition, that if $l_{\theta} \equiv log h_{\theta}$, $\big| \frac{\partial^2}{\partial \theta^2} l_\theta(x, y) \big| \leq K < \infty$ for all $|\theta|\leq \epsilon$, some $\epsilon>0, K <\infty$,

Then, (2.22) holds.
\end{proposition}

\textbf{Remark: } Using the techniques of Le Cam, the conditions may be relaxed to,

1) $h_{\theta}$ QMD with $\dot{l}_{\theta}(X,Y)=a(X)+b(Y)$

2) $\int exp[\theta(a(x)+b(y))]h_{\theta}(x,y)dxdy<\infty, |\theta|\leq \epsilon, \epsilon>0$

The techniques this relaxation requires may obscure the classical argument.

This means that the variational distance between $P_{\theta,n}$, the joint distribution of $(X_i,Y_i),i=1,...,n$ under $\mathcal{M}$, and $\tilde{P}_{\theta,n}$, the joint distribution of the $(X_i,Y_i),i=1,...,n$ under $\tilde{\mathcal{M}}$, goes to $0$ uniformly for $|\theta|\leq \frac{M}{\sqrt{n}}$, all $M<\infty$.

This implies that,
\begin{align*}
    P_{\theta_n n}[\sqrt{n}\hat{C}_n \geq t_n]=\tilde{P}_{\theta_n n}[\sqrt{n}\hat{C}_n \geq t_n]+o(1)
\end{align*}

for any sequence ${t_n}$, and $\theta_n=\frac{M}{\sqrt{n}}$.

Statistically, we can interpret this as follows. If $\hat{C}_n$ yields significant evidence against $H$ under $\mathcal{M}$, and $\tilde{\mathcal{M}}$ holds, we are being misled since, under $\tilde{\mathcal{M}}$, $X$ and $Y$ are always independent. But if we are convinced that $\mathcal{M}$ holds, we would not be using ranks.

In Theorem 2.2, we shall show that this conclusion holds generally for rank statistics which like $\hat{C}_n$
 are getting at the "local dependence" of $X$ and $Y$.
 
 These remarks can be viewed in another way in the context of semiparametric copula models.

\subsubsection{Semiparametric Copula models}\label{section:semiparametric}

 We start with a parametric model of densities $\mc{M}=\{h_\theta : |\theta| < 1\}$, where $h_\theta : (0, 1)^2 \ra [0, \infty)$ for all $|\theta| < 1$,  and $h_0(x, y) \equiv 1$. Let
\[ \mc{F} := \{a : [0, 1] \ra [0, 1] \text{ absolutely continuous}, a' > 0, a(0) = 0, a(1) = 1 \}.\]
i.e., $\mc{F}$ is the set of all absolutely continuous, strictly increasing transformations which map $0$ to $0$ and $1$ to $1$. Given $q, r \in \mc{F}$, we may define
\[ h_\theta(x, y, q, r) := h_\theta(q(x), r(y)) q'(x) r'(y), ~~ (x, y) \in (0, 1)^2. \]
Note that if $h_\theta$ is the density of $(X, Y)$, then $h_\theta(\cdot, \cdot, q, r)$ is the density of the pair $(q^{-1}(X), r^{-1}(Y))$. Define a semiparametric model by
\begin{align*}
   \mc{P} := \{h_\theta(x, y, q, r) : |\theta| < 1, q, r \in \mc{F} \}. 
\end{align*}

Note that we can, without loss of generality, take the domain of $h_\theta$ to be $R^2$ rather than $[0,1]^2$ since any statistical points we shall make are unchanged if we consider as basic alternatives $(F_1(X),F_2(Y))$ where $F_1,F_2$ are the cdf of $X,Y$ respectively, under $H$, rather than general $X,Y$. As in the past we shall call this vector, $(U,V)$.

Note also that we could similarly consider maps $a(.)$ belonging to $\mc{F}$ as having range $R$ rather than $(0,1)$. These models are natural candidates for using ranks in relation to any inference for $\theta$ since the "scales" a and b are unknown.

Note that $h_0(x, y, q, r) = q'(x) r'(y)$ corresponds to independence for all $(q, r)$ pairs.

Let for $F_{\theta},G_{\theta}$ the cdf of $X,Y$ respectively, 
\begin{align} \label{equation 2.23}
 \mc{Q} := \{Q_\theta : (F_\theta(X), G_\theta(Y)) \sim Q_\theta, ~(X, Y) \sim h_\theta, |\theta| < 1\}.    
\end{align}

Then $\mc{Q}$ is a parametric model which generates the same semiparametric copula model as $\{h_\theta(x,y) : |\theta| < 1\}$.

A parametric copula model is one where, for every $\theta$, the marginal distributions of $(X,Y)$ are $\mathcal{U}(0,1)$.

All classical parametric copula models, Archimedean, Gaussian, etc. are describable as in (2.23) with, in general, multivariate $\theta$.

As usual in semiparametric theory we now carry out calculations for parametric submodels of $\mc{P}$ where $q$ and $r$ depend on $\theta$ as well as $h_\theta(x,y)$. We will argue that $\mathcal{Q}$ is a least favorable submodel for estimating $\theta$ in $\mc{P}$.

Consider a parametric submodel, with $\bar{\theta} = (\theta_1, \theta_2, \theta_3)$, and
\[ \dens_{\bar{\theta}}(x, y) := h_{\theta_1}(x, y, q_{\theta_2}, r_{\theta_3}), ~~ |\bar{\theta}|_\infty < 1. \]
Let $\ell_{\bar{\theta}} := \log \dens_{\bar{\theta}}$, and let
\[\nabla \ell_{\bar{0}} (x, y) := (D_{1}(x, y), D_{2}(x, y), D_{3}(x, y)) \]
be the gradient of $\ell$ (in the $\bar{\theta}$ variable) at $\bar{\theta} = 0$. We suppose that the family $\{\dens_{\bar{\theta}}: |\bar{\theta}|_\infty < 1\}$ is quadratic mean differentiable at $\bar{\theta} = \bar{0}$. That is, it satisfies
\beq\label{eq:QMD}\tag{QMD} \E_{\bar{0}} \big(\dens_{\bar{\theta}}(X, Y)^{1/2} - \dens_{\bar{0}}(X, Y)^{1/2}- (1/2) \dens_{\bar{0}}(X, Y)^{1/2} (\nabla \ell_{\bar{0}}(X, Y), \bar{\theta})\big)^2 = o(|\bar{\theta}|_2^2),\eeq
where $\E_{\bar{0}}$ signifies that the random variables $(X, Y)$ are distributed according to $\dens_{\bar{0}}(x, y) = q_0'(x) r_0'(y)$. Note that
\[ D_{1}(X, Y) = \frac{\partial}{\partial \theta_1} \log h_{\theta_1}(q_{\theta_2}(X), r_{\theta_3}(Y)) \bigg|_{\bar{\theta}= 0} = \dot{\ell}_0(q_0(X), r_0(Y)), \]
(here $\dot{\ell}_0 := \ptl \log h_\theta / \ptl \theta |_{\theta = 0}$) and since $h_{\bar{0}} \equiv 1$, we have
\[ D_{2} (X, Y) = \frac{\ptl}{\ptl \theta_2} \log q_{\theta_2}'(X) \bigg|_{\theta_2 = 0}, ~~~ D_{3}(X, Y) = \frac{\ptl}{\ptl \theta_3} \log r_{\theta_3}'(Y)  \bigg|_{\theta_3 = 0}.\]

$\bold{Remarks: }$ 1) Let 
\beq\label{eq:q-theta-r-theta-def}
\begin{aligned}
q_\theta'(x) &:= q_0'(x) \big(1 + \theta a(q_0(x))\big), \\
r_\theta'(y) &:= r_0'(y) \big(1 + \theta b(r_0(y))\big), 
\end{aligned}
\eeq
where $|a| < 1, |b| < 1$, and $\E a(U) = \E b(U) = 0$, for $U \sim \mathrm{Unif}(0, 1)$. If $\{h_\theta : |\theta| < 1\}$ satisfies assumptions I and II, then the parametric copula model defined by (2.23) above satisfies \eqref{eq:QMD}.

The tangent spaces (using the notation of \cite{BickelEAES1998}) of the model defined by \eqref{eq:q-theta-r-theta-def} at $\theta_1 = \theta_2 = \theta_3 = 0$ are
\[ \dot{\mc{P}}_{\theta_1} = [\dot{\ell}_0(q_0(X), r_0(Y))],~~ \dot{\mc{P}}_{\theta_2} = [a(q_0(X))], ~~ \dot{\mc{P}}_{\theta_3} = [b(r_0(Y))]. \]
Since the set of $a, b$ as above are dense in $L_2^0([0, 1]) := \{c(U) : \E [c(U)^2] < \infty, \E c(U) = 0\}$, the tangent spaces of the semiparametric copula model are

\begin{align*}\label{equation 2.24}
   \dot{\mc{P}}_{\theta} = [\dot{\ell}_0(q_0(X), r_0(Y))],
\end{align*}
\begin{align*}
   \dot{\mc{P}}_q = \{f(q_0(X)) : f \in L_2^0([0, 1])\}, ~~ \dot{\mc{P}}_r = \{f(r_0(Y)) : f \in L_2^0([0, 1])\}.  
\end{align*}

While $\dot{\ell}_0(q_0(X), r_0(Y))$ is the score function for the base model $\{h_\theta(\cdot, \cdot, q_0, r_0) : |\theta| < 1\}$, we now obtain that
\begin{align}
    \dot{\ell}_0^*(X, Y) := \dot{\ell}_0(q_0(X), r_0(Y)) - \E[\dot{\ell}_0(q_0(X), r_0(Y)) ~|~ X] - \E[\dot{\ell}_0(q_0(X), r_0(Y)) ~|~ Y]
\end{align} 
is the efficient score function for the model $\mc{P}$, since $\E[\dot{\ell}_0(q_0(X), r_0(Y)) ~|~ X]$ is the projection of $\dot{\ell}_0(q_0(X), r_0(Y))$ on $\dot{\mc{P}}_q$ and similarly for $\dot{\mc{P}}_r$ (see \cite[Chapter 3]{BickelEAES1998}). 

Consider the submodel of $\mc{P}$, given by \eqref{eq:q-theta-r-theta-def}, with $\theta_1 = \theta_2 = \theta_3 = \theta$, and
\[ a(q_0(X)) = -\E [\dot{\ell}_0(q_0(X), r_0(Y)) ~|~ X], ~~ b(r_0(Y)) = -\E[\dot{\ell}_0(q_0(X), r_0(Y)) ~|~ Y]. \]
Given assumptions I and II, we have that \eqref{eq:QMD} holds, and hence this model is also quadratic mean differentiable at $\theta = 0$, with score function $\dot{\ell}_0^*$. For estimation, this means that if an estimate is regular and linear on $\mc{P}$, then its asymptotic variance $\sigma^2 / n \geq (I^*)^{-1} / n$, where
\[ I^* := \E_0[\dot{\ell}_0^*(q_0(X), r_0(Y))^2]. \]
For a rank test statistic $T_n$ such that under $\theta = 0$, we have that
\begin{align*}\label{eq:semiparametric-joint-clt} 
\bigg(\sqrt{n}T_n, \frac{1}{\sqrt{n}} \sum_{i=1}^n  \dot{\ell}_0^*(X_i, Y_i) \bigg) \stackrel{d}{\lra} N\bigg(\begin{pmatrix} 0\\ 0 \end{pmatrix},  \begin{pmatrix} \sigma^2 & \rho^* \sigma \sqrt{I^*} \\ \rho^* \sigma \sqrt{I^*} & I^* \end{pmatrix}\bigg), 
\end{align*}
Then as in (2.17), the Pitman efficiency of $T_n$ with respect to the optimal asymptotic test is

\begin{align*}
    e^*(T_n) = (\rho^*)^2.
\end{align*}

2) We can exhibit explicit families $\{q_\theta : |\theta| < 1\}$, $\{r_\theta : |\theta| < 1\}$ such that the parametric model $\{h_\theta(\cdot, \cdot, q_\theta, r_\theta) : |\theta| < 1\}$ has score function $\dot{\ell}_0^*$. Set $q_\theta := F_\theta^{-1}$, $r_\theta := G_\theta^{-1}$, where $F_\theta, G_\theta$ are the respective marginal cdfs of $X, Y$ under $h_\theta$. (Note then the distribution $h_\theta(\cdot, \cdot, q_\theta, r_\theta)$ is given by the law of $(F_\theta(X), G_\theta(Y))$, if $(X, Y)$ is distributed according to $h_\theta$. So we are essentially going from the semiparametric copula generated model back to the parametric copula model.) Let $\dot{s}_0$ denote the score function at $\theta = 0$ for this model. A change of variables calculation then gives that 
\[ \dot{s}_0(X, Y) = \dot{\ell}_0(q_0(X), r_0(Y)) - \frac{1}{f_0(X)} \frac{\ptl f_\theta}{\ptl \theta} (q_0(X)) \bigg|_{\theta = 0} - \frac{1}{g_0(Y)} \frac{\ptl g_\theta}{\ptl \theta}(r_0(Y)) \bigg|_{\theta = 0}. \]

From semi parametric theory it follows that the efficient score function at $\theta=0$ is the projection of $\dot{\ell}_0$ on the ortho complement of the nuisance tangent space $\dot{\mathcal{P}}_q+\dot{\mathcal{P}}_r$. Without loss of generality, we can take $q_0$ and $r_0$ the identity. Thus $\dot{\mathcal{P}}_q+\dot{\mathcal{P}}_r=\Lambda_1$ as defined in (2.13) for n=1.

\subsubsection{ Local power behaviour of classical statistics}

The properties of $\hat{T}_{BKR}$, $\hat{T}_{S}$ are essentially well known. 

It is obvious, for instance, from the expansion

\begin{align}
    \hat{T}_{S}=\frac{12}{n}\sum_{i=1}^{n}{(U_i-\frac{1}{2})(V_i-\frac{1}{2})}+o_p(n^{-\frac{1}{2}}).
\end{align}

that $\hat{T}_{S}$ has positive Pitman efficiency in any $\mc{M}$, unless $\dot{l}_{0}(X,Y) \perp (U-\frac{1}{2})(V-\frac{1}{2})$.

Nevertheless, the following theorem may not be known in full generality.

\begin{theorem}
If $\dot{l}_0(X,Y) \neq a(X)+b(Y)$ for some $a,b$, then $\hat{T}_{BKR}$ is asymptotically unbiased.
\end{theorem}

As we have argued in the previous section, given a semiparametric formulation of a model, a score function in $\Lambda_1$,  corresponds to models in which dependence cannot be measured reliably.

$\bold{Proof: }$ By general theory, under $H$,

\begin{align*}
    n\hat{T}_{BKR} \stackrel{d}{\lra} \int [T^{*}]^2(\tend{u}) d\tend{u}.
\end{align*}

and, under $\{P_{n\theta_n}\}$, if $(T_n(.),\frac{1}{\sqrt{n}}\sum_{i=1}^n l_0^{*}(X_i,Y_i)) \stackrel{d}{\lra} (T^{*}(.),Z)$,

\begin{align*}
    n\hat{T}_{BKR} \stackrel{d}{\lra}  \int (T^*(\tend{u})+tC(\tend{u}))^2 d\tend{u}.
\end{align*}

where 

\begin{align*}
    C(\tend{u})=cov(T^*(\tend{u}),Z)
\end{align*}

By classical theory, we can represent, for $Z_{ij}$ iid $N(0,1)$,

\begin{align*}
    \int [T^{*}]^2(\tend{u}) d\tend{u} &= \sum_{i,j\geq 1} \lambda_{ij} Z_{ij}^2
 \end{align*} 
 
 where 
 
 \begin{align}
 \sqrt{ \lambda_{ij}}Z_{ij}=\int T^*(\tend{u}) f_{ij}(\tend{u}) d \tend{u}
 \end{align}
 
 and, 
 
 \begin{align*}
     f_{ij}(\tend{u})=2sin(\pi i u)sin(\pi j v).
 \end{align*}
 
 Then,
 
\begin{align*}
    \int (T^*(\tend{u})+tC(\tend{u}))^2 d\tend{u} &\sim \sum_{i,j\geq 1} \lambda_{ij} (Z_{ij}+tC_{ij})^2.
\end{align*}

Here, as shown by Blum et al (1961) \cite{BlumDFTIJ1961}, $\lambda_{ij}=\frac{1}{\pi^4 i^2 j^2}$,

\begin{align*}
    C_{ij}=\int C(\tend{u}) f_{ij}(\tend{u}) d\tend{u}
\end{align*}

and $f_{ij}$ are the eigenfunctions of the covariance kernel given in (2.14). By Lemma A1, $C(.) \not\equiv 0$. By Zygmund \cite{ZygmundTS2002} Problems 1.8.4 and 1.8.6, the $f_{ij}$ are a complete orthonormal system in $L_2[0,1]^2$. Therefore, if $C(.) \not\equiv 0$, some $C_{ij} \neq 0$. Then, if $t\neq 0$, the distribution of $(Z_{ij}+tC_{ij})^2$ is stochastically larger than that of $Z_{ij}^2$ and the theorem follows. 
\qedsymbol

As, is well known in the context of goodness of fit tests \cite{LehmannTSHM2006a}, no test sequence can have substantial local power in more than a finite set of directions, essentially choices of $\dot{l}_0$ in our context. This carries over to semiparametric hypotheses within nonparametric models, such as independence generally, see for instance Bickel, Ritov and Stoker (2006)\cite{BickelTTGFA2006}. This is clear from consideration of $\hat{T}_{BKR}$, where it is natural to take as primitive directions, the elements of the Fourier basis, that is, $\dot{l}_0 (F_1^{-1}(U),F_2^{-1}(V))=f_{ij}(U,V)$. The weight $\frac{1}{\pi^4 i^2 j^2}$ placed on each of these in the limit distribution of $\hat{T}_{BKR}$ diminishes as a function of $i,j$, corresponding to their "wiggliness". Bickel, Ritov and Stoker (2006) \cite{BickelTTGFA2006} discuss various ways of favoring preferred directions while retaining  some  local power in all directions. We will turn to this topic in the contex of the $f_{ij}$ and other orthonormal bases elsewhere.

In the remainder of this section, we will

1) Make explicit the construction of rank statistics corresponding to the same statistics in terms of the $(U_i,V_i)$ (which are not observed)

2)Indicate how such statistics can be used for a simpler purpose, approximating $\hat{T}_{BKR}$ to arbitrary precision in $O(n)$ operations, as well as building blocks for selective test of independence.

We use the following lemma, proved in the Appendix.

\begin{lemma}
Let $a,b \in L_2[0,1]$. Let $A(u)\equiv \int_0^u a(s)ds$, $B(v)\equiv \int_0^v b(t)dt$. Let $(R_i,S_i),i=1,...,n$ be the ranks of the $(U_i,V_i)$ coordinatewise where $U_i=F_1(X_i), V_i=F_2(Y_i)$. Let $\bar{A}(u)=A(u)-EA(u)$, $\bar{B}(v)=B(v)-EB(v)$. Then, under $H$,

\begin{align}
    \int_0^1\int_0^1 T_n(u,v)a(u)b(v)dudv=\frac{1}{\sqrt{n}}\sum_{k=1}^n{\bar{A}(\frac{R_k}{n})\bar{B}(\frac{S_k}{n})}+o_p(n^{-\frac{1}{2}})
\end{align}
\end{lemma}

$\bold{Remark: }$ If $a$ is bounded, $EA(U)=\frac{1}{n}\sum_{k=1}^n{A(\frac{k}{n})}+O(\frac{1}{n})$

Now apply Lemma 2.2 to the $f_{ij}$ and obtain,

\begin{align}
    \hat{T}_{ij}\equiv \int T_n(\tend{u})f_{ij}(\tend{u})d\tend{u}=\frac{2}{\sqrt{n}}\sum_{k=1}^n\bigg[\frac{cos(\frac{\pi i R_k}{n})cos(\frac{\pi j S_k}{n})}{(\pi ij)^2}\bigg]+o_p(n^{-\frac{1}{2}})
\end{align}

since $\int_0^u sin(\pi i t)dt=\frac{1-cos(\pi iu)}{\pi i}$ and $\int_0^1 cos(\pi i u)du=0$.

By Proposition 2.1, $\{\hat{T}_{ij}: 1\leq i,j \leq N\}$ are asymptotically jointly Gaussian with mean 0, independent, and with variance $\{[\pi^4i^2j^2]^{-1}:i,j\leq N\}$. The independence holds because $f_{ij}(.)$ are an orthornormal basis of $L_2[0,1]^2$, eigenfunctions of the covariance kernel of $T^*$ and (2.27).

Moreover, by standard Fourier analysis,

\begin{align*}
    T_n(\tend{u})=\sum_{i,j\geq 1}{\hat{T}_{ij}f_{ij}(\tend{u})} \textit{ a.e }
\end{align*}

and

\begin{align}
    n\hat{T}_{BKR}=\int T_n^2(\tend{u})d\tend{u}=\sum_{i,j\geq 1} \hat{T}_{ij}^2
\end{align}

If
\begin{align}
    \hat{T}_{n,M}^2(\tend{u})\equiv \sum_{i,j= 1}^M \hat{T}_{ij} f_{ij}(\tend{u})
\end{align}

then,

\begin{align*}
    \int T_{n,M}^2(\tend{u})d\tend{u}=\sum_{i,j= 1}^{M}\hat{T}_{ij}^2
\end{align*}

and

\begin{align*}
    \int(T_{n,M}(\tend{u})-T_n(\tend{u}))^2 d\tend{u} \to 0 \textit{ in probability as } M\to \infty
\end{align*}

We have proved,

$\bold{Proposition 2.5: }$ Under contiguous alternatives, $n\hat{T}_{BKR}$ can be approximated arbitrarily closely by $
\sum_{i,j\geq 1}^M \hat{T}_{ij}^2$ and these statistics can be computed in $O(n)$ steps.

In fact we can use the $\hat{T}_{ij}$ for selective inference since they are symptotically independent. We can use them in lexicographic order in the $(i,j)$ basis, corresponding to their wiggliness. Note again that we are using the fact that $sin(\pi i u)sin(\pi j v)$ are both eigenfunctions of the covariance kernel of $T_n(.,.)$ and an orthonormal basis of $L_2[0,1]^2$.

It is attractive to consider alternative orthonormal basis which localize failure of independence more specifically. A natural candidate based on Rademacher functions, is given as follows. Choose a scale of interest, $\Delta=2^{-N}$. Let

\begin{align}
\Psi_p(u)=\Delta\Psi(\frac{u}{\Delta}-p), \quad p=0,...,2^N-1    
\end{align}

where

\begin{align*}
    \Psi(u)=sgn(u-\frac{1}{2})
\end{align*}

By construction, the $\Psi_p(.)$ are orthonormal. We could take $\Delta >0$ generally, but our choice enables us to vary $N=0,1,2,...$ and create a full orthonormal basis of 

\begin{align*}
    L_2^0(0,1) \equiv \{h\in L_2(0,1): \int h(u)du=0\}.
\end{align*}

Let,

\begin{align*}
    \hat{S}_{N,p_1,p_2}&=\int_0^1\int_0^1 T_n(u,v) \Psi_{p_1}(u) \Psi_{p_2}(v) dudv\\
    &=\frac{1}{n}\sum_{i=1}^n \lambda_{p_1}(\frac{R_i}{n}) \lambda_{p_2}(\frac{S_i}{n})+o_p(n^{-\frac{1}{2}})
\end{align*}

under $H$, where

\begin{align*}
    \lambda_p(u)&=2^{-N}\lambda_0(2^N u-p)\\
    \lambda_0(u)&=\int_u^1 \Psi(s)ds-\frac{1}{4}\\
    &=t1(0\leq t\leq \frac{1}{2})+(1-t)1(\frac{1}{2}\leq t\leq 1)-\frac{1}{4}
\end{align*}
Under $H$, the $\hat{S}_{N,p_1,p_2}$ are asymptotically Gaussian with mean 0, (jointly), iid for fixed N with variance $\frac{2^{-N}}{48}$, but dependent as N varies.

Selective inference is   done for successive scales. The approach can be naturally  extended to independence of more than 2 variables. An approach related to this proposal has been put forward by Zhang \cite{ZhangBIO2019}.

\subsection{Local independence and a family of test statistics and parameters related to $\hat{C}_n$}

Chatterjee's $\hat{C_n}$ and recent work of Lin and Han \cite{LinBPCRA2021} and Deb, Ghosal and Sen \cite{DebMATSO2020} suggest a family of parameters based on a representation of the joint distribution of a pair $(X,Y)$ with $F$ continuous which leads to novel parameters measuring dependence and test statistics for H.

If F is the probability distributions of a pair  $(X,Y)$, we can consider it in the form, $(F_1(.),F_2(.|.))$ where $F_2(y|x)$ is the cdf of the regular conditional probability distribution of $Y$ given $X$. We can extend this representation by considering $Y^{(1)},...,Y^{(m)}$ conditionally iid with distribution $F_2(.|X)$. Call the joint distribution of $(X,Y^{(1)},...,Y^{(m)})$, $F^{(m)}$. With this notation, the Chatterjee and Dette et al. parameter in form (1.2) if F is continuous is,

\begin{align}
    T(F)= 1-3 E_{F^(2)}\{ |F_2(Y^{(1)})-F_2(Y^{(2)})| \}
\end{align}

This point of view leads naturally to considering the process defined on $[0,1]^3$ by,

\begin{align*}
    \hat{C_L}(u,v_1,v_2) \equiv \frac{1}{n-1)} \sum_{i=1}^{n-1} 1(\hat{F_1}(X_{(i)}) \leq u, \hat{F_2}(Y_{[i]}) \leq v_1, \hat{F_2}(Y_{[i+1]}) \leq v_2  )
\end{align*}
where, if $(X_1,Y_1),...,(X_n,Y_n)$ are the original sample, $X_{(1)} < ... <X_{(n)}$ are the ordered $X_i$ and $(X_{(i)},Y_{[i]}) = (X_{D_i}, Y_{D_i})$ where $D_i=j$ iff $X_{(i)}=X_j$, and $\hat{F}_j$ are the usual empirical marginal cdf.

Thus, in the notation of (1.1),
\begin{align*}
    \hat{C_L}(u,v_1,v_2) =\frac{1}{n-1} \sum_{i=1}^{[nu]} 1(\frac{r_i}{n} \leq v_1, \frac{r_{i+1}}{n} \leq v_2).
\end{align*}

Next, define in analogy to $T_n(.,.)$,

\begin{align*}
    S_n(u,v_1,v_2) = \sqrt{n} (\hat{C_L}(u,v_1,v_2)- C_L(u,v_1,v_2))
\end{align*}

where,

\begin{align*}
    C_L(u,v_1,v_2)=\int_{0}^{u} F_{2}^{*}(v_1|Q_1(s)) F_{2}^{*}(v_2|Q_1(s)) ds
\end{align*}

and $F_{2}^{*}$ is the conditional cdf of $F_{2}(Y)$ given $X=x$.

In the case we pursue, $X,Y$ independent and $F$ continuous, $C_L(u,v_1,v_2)=u v_1 v_2$ and we focus on,

\begin{align}
    S_n(v_1,v_2)=\sqrt{n} (\hat{C_L}(v_1,v_2)-v_1v_2)
\end{align}

where $\hat{C_L}(v_1,v_2) \equiv \hat{C_L}(1,v_1,v_2)$.

Rank statistics such as $\hat{T}_{BKR}, \hat{T}_{S}$ have analogues, (neglecting terms of $o(\frac{1}{n})$).

In particular, define

\begin{align}
    \hat{T}_{S}^{(L)}=12(\int_{0}^{1} \int_{0}^{1} v_1v_2d\hat{C_L}(v_1,v_2) -\frac{1}{4}) \notag\\
    =\frac{12}{n-1} \sum_{i=1}^{n-1}(\frac{r_i}{n}\frac{r_{i+1}}{n}-\frac{1}{4})
\end{align}

All such statistics can evidently be computed using the $(r_i,r_{i+1})$ just as their global analogues use $(R_i,S_i)$.

$\bold{Remark:}$ Such processes and  statistics have been defined previously in a different context, see for example, Genest et al (2007)\cite{GenestREBDD2007}.

We shall see in Theorem 2.2 that $\hat{T}_{S}^{(L)}$ and all statistics based on $\hat{C}_L$ have the same poor local power behaviour we claimed for $\hat{C}_n$.

The following result will prove useful. Its proof is in the appendix.

\begin{lemma}
Suppose $(U_i,V_i),i=1,...,n$ are i.i.d $\mathcal{U}(0,1)$. Define, as below, $V_{[i]}=V_j$ iff $U_j=U_{(i)}$ where $U_{(1)} <...< U_{(n)}$ are the ordered $U_i$. Given $g\in L_2[0,1]^2$, define $g_1(s)\equiv \int_0^1 g(s,t)dt$, \newline
$g_2(t)\equiv \int_0^1 g(s,t)ds$. Let $a,h^{(1)},h^{(2)} \in L_2[0,1]^2$ and use the subscript notation defined, noting that $h_1^{(j)}(V_1)=E[h^{(j)}(V_1,V_2)|V_1], h_2^{(j)}(V_1)=E[h^{(j)}(V_2,V_1)|V_1]$, $j=1,2$ etc.

Let,

\begin{align*}
    A_n \equiv \frac{1}{\sqrt{n}}\sum_{k=1}^n a(U_{(k)},V_{[k]})\\
    H_n^{(i)} \equiv \frac{1}{\sqrt{n}}\sum_{k=1}^{n-1} h^{(i)}(V_{[k]},V_{[k+1]})
\end{align*}

Then, $(A_n, H_n^{(1)}, H_n^{(2)}) \stackrel{d}{\lra} (Z_0,Z_1,Z_2) $

where $(Z_0,Z_1,Z_2)$ are jointly Gaussian, mean 0 and

\begin{align*}
    Var(Z_0)=Var (a(U,V))
\end{align*}   
 
 \begin{align}
    &Var(Z_1)=Var(h^{(1)}(V_1,V_2))+Var(h_1^{(1)}(V_1))+Var(h_2^{(1)}(V_2)) \\
    &Cov(Z_0,Z_j)=Cov(a_2(V_1),h_1^{(j)}(V_1)+h_2^{(j)}(V_1)),\quad j=1,2.\\
    &Cov(Z_1,Z_2)=Cov(h^{(1)}(V_1,V_2),h^{(2)}(V_1,V_2))+Cov(h_1^{(1)}(V_1),h_1^{(2)}(V_1))+Cov(h_2^{(1)}(V_2),h_2^{(2)}(V_2))
 \end{align}   

\end{lemma}

Let 

\begin{align*}
h(\tend{v},w_1,w_2)=1(w_1\leq v_1)1(w_2\leq v_2)-v_1v_2    
\end{align*}
 
and define, consistent with our general notation,

\begin{align*}
    h_1(\tend{v},w_1)=\int_0^1 h(\tend{v},w_1,w_2)dw_2=(1(w_1\leq v_1)-v_1)v_2
\end{align*}

and $h_2(\tend{v},w_2)$ analogously,

Then,

\begin{align*}
    S_n(v_1,v_2)=\frac{1}{\sqrt{n}} \sum_{k=1}^{n-1} h(\tend{v},\frac{r_k}{n},\frac{r_{k+1}}{n})
\end{align*}

Let 
\begin{align}
    \tilde{h}(\tend{v},w_1,w_2)= h(\tend{v},w_1,w_2) - h_1(\tend{v},w_1)-h_2(\tend{v},w_2)
\end{align}

Define,

\begin{align*}
    R_n(\tend{v}) &\equiv \frac{1}{\sqrt{n}} \sum_{i=1}^{n-1} h(\tend{v},V_{[i]},V_{[i+1]}) \\
    Q_n(\tend{v}) &\equiv \frac{1}{\sqrt{n}} \sum_{i=1}^{n-1} \tilde{h}(\tend{v},V_{[i]},V_{[i+1]})\\
    \Lambda_n^*(\tend{v}) &\equiv \frac{1}{\sqrt{n}} \sum_{i=1}^{n-1} (h(\tend{v},V_{[i]},V_{[i+1]})-\tilde{h}(\tend{v},V_{[i]},V_{[i+1]}))\\
    &=\frac{1}{\sqrt{n}}\sum_{i=1}^{n-1} (h_1(\tend{v},V_{[i]}))+h_2(\tend{v},V_{[i+1]}))
\end{align*}

Also, by definition,

\begin{align*}
    Q_n(\tend{v})=R_n(\tend{v})-\Lambda_n^*(\tend{v}).
\end{align*}

Here is an analogue to Proposition 2.1.

\begin{theorem}
If $X$ and $Y$ are independent and $F_1, F_2$ are continuous, then,

\begin{enumerate}[label=(\roman*)]
    \item $S_n(v_1,v_2)=Q_n(v_1,v_2)+o_p(1)$\\
    \item $\Lambda_n^* \stackrel{d}{\lra} \Lambda^*$
    \\
    and
    \item $S_n(.,.) \stackrel{d}{\lra} S(.,.)$
\end{enumerate}
\end{theorem}

where $S(.,.)$ is a Gaussian field with mean 0 and

if $\tend{a}=(a_1,a_2), \tend{b}=(b_1,b_2)$, then

\begin{align}
    Cov(S(\tend{a}),S(\tend{b})) =Cov(T^*(\tend{a}),T^*(\tend{b}))+Cov(\Lambda^*(\tend{a}),\Lambda^*(\tend{b}))+o(1)
\end{align}

where $T^*$ is given in (2.14), and $\Lambda^*(v_1,v_2) \sim v_1 E_2(v_2)+v_2E_2(v_1)$.

$\bold{Remark: }$ We can interpret (2.35) as,

\begin{align*}
    S(\tend{a}) \sim T^*(\tend{a}) +  \Lambda^*(\tend{a}) 
\end{align*}

where $T^*$ and $\Lambda$, are $\bold{independent}$ mean 0 Gaussian processes and $T^*(.)$ has the covariance given in (2.14). Here "$\sim$" indicate distributional equivalence since the expression arises from the CLT for 2 dependent variables Orey(1958) \cite{OreyCLTMD1958}.

$\bold{Proof: }$ As with Theorem 2.1, we need to show tightness and joint FIDI convergence of the 2 processes,$S_{n}-\Lambda^*_{n}$ and $\Lambda^*_{n}$, and prove they have continuous sample surfaces. The delta method is applied as in Theorem 2.1.

Given $\{U_{(1)} <...<U_{(n)}\}\equiv \mathcal{B}_n$, the $V_{[i]}$ are iid $\mathcal{U}(0,1)$ and we can just write $V_i$.

Apply the 2 dependent CLT (Orey(1958)\cite{OreyCLTMD1958}) to obtain that, since \newline $E[Q_n(.)|\mathcal{B}_n]=E[\Lambda_n^*(.)|\mathcal{B}_n]=0$,

\begin{align*}
    (Q_n(.), \Lambda_n^*(.)) \to (Q,\Lambda^*)
\end{align*}

where $Q$ and $\Lambda^*$ are independent, and the covariances are as specified by applying Lemma 2.1.

$\{\Lambda_n^*\}$ is evidently tight. $Q_n(.)$ is also by a direct calculation showing that

\begin{align}
    E(\Delta_{\epsilon} Q_n)^2(\tend{u}) (\Delta_{\epsilon} Q_n)^2(\tend{u}+\tend{\epsilon}) \leq C \epsilon^2
\end{align}

where $\tend{\epsilon}=(\epsilon_1,\epsilon_2)$ and 

\begin{align*}
    \Delta_{\epsilon}G(u_1,u_2) \equiv G(u_1+\epsilon_1,u_2+\epsilon_2)-G(u_1+\epsilon_1,u_2)-G(u_1,u_2+\epsilon_2)+G(u_1,u_2) \qquad \qedsymbol
\end{align*}

We generalize Proposition 2.3.

\begin{theorem}
Suppose $H$ holds and the condition of Proposition 2.3 is satisfied. 

Then, $S_n(.)$ is asymptotically independent of $L_n \equiv \frac{1}{\sqrt{n}} \sum_{i=1}^n \dot{l}_0(X_i,Y_i)$. Thus, any test statistic based on $S_n(.)$ which has a joint limit with $L_n$ has no power.
\end{theorem}

\textbf{Proof: } In view of Theorem 2.2, we need only check that,
\begin{align}
    \lim_n cov_0 \bigg( \frac{1}{\sqrt{n}}\sum_{i=1}^n\big(h_1(\tend{v},V_{[i]})+h_2(\tend{v},V_{[i+1]})\big),L_n \bigg)=0
\end{align}

and

\begin{align}
    \lim_n cov_0 \bigg( \frac{1}{\sqrt{n}}\sum_{i=1}^n \tilde{h}(\tend{v},V_{[i]},V_{[i+1]}),L_n \bigg)=0
\end{align}

since joint Gaussianity of these statistics follows from Orey's theorem once we write (with $X=U,Y=V$)

\begin{align*}
    L_n=\frac{1}{\sqrt{n}}\sum_{i=1}^n \dot{l}_0(U_{[i]},V_{[i]})
\end{align*}

(2.42) follows by substituting $\dot{l}_0$ for $a(.,.)$ and $h_j$ for $h^{(j)}$ in Lemma 2.3.

By assumption, $a_2=0$. Similarly, by substituting $\tilde{h}$ by $h_j$, (2.43) follows.\qedsymbol

\subsection{Discussion}

1) If F is not continuous, we can have $(X_1,Y_1), (X_2,Y_2)$ such that $X_1=X_2$ or $Y_1=Y_2$ or both. In those cases, the distribution of global rank statistics under H depends on the underlying $F_1, F_2$. For local statistics, as Chatterjee points out, ties occurring only in X can be broken randomly.

Valid tests can be carried out by conditioning on $\hat{F}_1(.)$ and $\hat{F}_2(.)$. For instance, if $(X,Y)$ is discrete, we would be as usual conditioning on the row sums of a contingency table.

2) Even in the continuous case, tests based on statistics such as $\hat{T}_{BKR}$ or $\hat{T}_{BKR}^{(L)}$are naturally carried out using Monte Carlo approximations to the null distributions using simulations of the $(U_i,V_i)$ under $H$.

\section{Functional Dependence}

In this section we will discuss property $C)$, functional dependence of the Chatterjee-Dette et al. parameter and a principle leading to $R,C$, and a number of other parameters defining functional dependence in different contexts. We change notation slightly here letting $F_X$ and $G_Y$ denote the marginal cdf of $X$ and $Y$.

\subsection{A Measure of functional dependence}

We begin with the solution to a simple problem. For $f_0(Y) \in L^2(Y)$, let 
\begin{align*}
    M(X,f_0(Y)):= sup\{ \rho(f_0(Y),g(X))^2,:g(X) \in L^2(X)\},
\end{align*}

Here $\rho$ is the correlation coefficient.

By the correlation inequality, 

\begin{align}
     M(X,f_0(Y))=\frac{Var(E(f_0(Y)|X))}{Var(f_0(Y))}.
\end{align}

This measure has the property that $M(X,f_0(Y))=1$ iff $f_0(Y)=h(X)$ a.s. for some $h$.

\subsection{C in this context}
  
The Chatterjee-Dette et al parameter emerges from this point of view as follows.
  
Note that for continuous $F_X,G_Y$, $Y=h(X)$ iff $G_Y(Y)=\tilde{h}(F_X(X))$ a.e. where

\begin{align*}
    \tilde{h}(x)=G_{Y}^{-1}\tilde{h}F_X(x)
\end{align*}

and $F^{-1}(u) \equiv inf\{x:F(x)\geq u\}$.

We obtain a measure of functional dependence of $Y$ on $X$ in (3.1) by taking $f_0(y) \equiv y$. Then,

\begin{align}
    M(X,Y) &= \frac{ Var(E(Y|X))}{ Var(Y)}\\
    &=1-\frac{E(Var(Y|X))}{Var(Y)}. \notag
\end{align}

This satisfies the functional dependence requirements C) but not B). Chatterjee and Dette et al. implicitly use that for $Y\in R$, independence is equivalent to

\begin{align}
    P(Y \leq y |X)=P(Y \leq y) \textit{ a.s.}
\end{align}
 
which corresponds to ,

\begin{align}
    M(X,1(Y \leq y))=\frac{E(Var(1(Y \leq y)|X))}{Var(1(Y \leq y))}=1
\end{align}

for almost all $y$.

Chatterjee then  achieved the goal of of characterizing independence by integrating the numerator and denominator  in (3.4) separately in y and defining C as in (1.2).

\subsection{Extensions}

A natural extension of (3.1) is to specify parameters having the extremal property for $g$ belonging to subsets of $L_2(X)$, and also for $Y,X$ with values in $R^p$ or more general spaces. A natural extension of this type is to define for $Y \in R^p$,

\begin{align}
M_L(X,Y) \equiv max\{ M(X,f(Y)),f \in \mc{F} \}    
\end{align}

where $\mc{F}=\{\textit{All linear functions of }Y \}$. It readily follows that,

\begin{align}
    M_L(X,Y)=\textit{Maximum eigenvalue of } \Sigma
\end{align}
where,

\begin{align*}
    \Sigma&=Var^{-\frac{1}{2}}(Y) Var(E(Y|X)) Var^{-\frac{1}{2}}(Y)\\
    &=I-Var^{-\frac{1}{2}}(Y) E(Var(Y|X)) Var^{-\frac{1}{2}}(Y)
\end{align*}

and $Var(W)$ is the variance-covariance matrix of $W \in R^p$.

This notion can be extended to $Y$ not fully dimensional by considering singular values and $Y$ not in $L_2$ by, as for $C$, considering $(G_{Y^{(1)}}(Y^{(1)}),...,G_{Y^{(p)}}(Y^{(p)}))$ where $Y\equiv (Y^{(1)},...,Y^{(p)})$.

The definition becomes valid even if $G_{Y^{(j)}}(.)$ are not continuous.

$\bold{\textit{The discrete case}}$

If $(X,Y)$ has finite support $x\in \{ x_1,...x_r\}, y\in \{ y_1,...y_s\}$, then, $M_L$ defined by (3.5) is the maximal canonical correlation between
\begin{align*}
(1(X=x_1),.....,1(X=x_r)) \textit{ and } (1(Y=y_1),.....,1(Y=y_s)).   
\end{align*}
In particular, $M_L$ and $R$ coincide, since any function of $(1(X=x_1),.....,1(X=x_r))$ can be written as a linear combination and similarly for $Y$. Using (1.3) we get,

\begin{align}
    M_L = \textit{maximal eigenvalue of } {\Vert \sum_{j=1}^{s} {p_{aj}q_{jb}} \Vert_{r\times r}}
\end{align}

(or the corresponding $s \times s$ matrix for $Y$).

Here, $p_{aj} \equiv P[Y=y_j|X=x_a], q_{jb} \equiv P[X=x_b | Y=y_j]$

There is an extensive literature on this measure and related topics in the contingency table literature, see Goodman(1965), Gilula and Haberman \cite{GilulaPFCPA1995} for instance.
\\
\\
\\
$\bold{\textit{Shape constrained subspaces of $L_2(X)$}}$ 
A natural extension to consider is to situations where we require that the dependence parameter $M$ should be 1 if and only if $Y=h(X)$ for $h \in \mc{H}$, a convex cone in $L_2(X)$. The simplest example, after $L_2$ itself which leads to $C$, is, the isotonic case, taking $\mc{H}=\{ \textit{Monotone functions of }X\}$. In fact, a class of such parameters exist, namely $T_{S}$ and related ones such as the normal scores correlation. All of these rely on the fact that, if $F_X$ and $G_Y$ are continuous, $Y=h(X)$ for $h$ monotone iff $F_X(X)=G_Y(Y)$ and the ordinary linear condition suffices. If we apply (3.1) and define,
\begin{align*}
    C_{\mc{H}}(X,Y) \equiv max\{\rho^2(Y,h(X)):h\in \mc{H}\},
\end{align*}

then, by [14] and (3.1),

\begin{align}
     C_{\mc{H}}=\frac{Var(\Pi_{\mc{H}_X}(Y))}{Var(Y)}
\end{align}

where,
\begin{align*}
    \Pi_{\mc{H}_X}(Y) \equiv \argmin\{E(Y-h(X))^2: h\in \mc{H}\}
\end{align*}
exists uniquely, if $\mc{H}$ is a closed convex cone in $L_2(X)$.

Cao and Bickel show, in \cite{CaoCTEPO2020}, that, in the isotonic case, unlike $C$, $C_{\mc{H}}$ can be estimated by plug in and the resulting $\hat{C}_{\mc{H}}$ is consistent and has interesting probabilistic features. Unfortunately its behaviour as a test of independence is worse than $\hat{T}_{S}$ which has the same functional dependence properties that $\hat{C}_{\mc{H}}$ has. 

For other $\mc{H}$, such as convex functions, plug in estimation seems impossible, see [13], as it is for $C$, and as with Dette et al's statistic is likely to be consistent only under strong conditions on $\Pi_{\mc{H}_X}(Y)$ as a function of $X$. We do not pursue these applications further.

\subsection{Estimation}

The passage from (1.2) back to (1.1) is obscure. Here is a possible clarification.

Chatterjee's approach relies on the second formula in (1.2). Specifically, he estimates $Var(1(Y\leq y)|X)$ by noting that, given $Y'=y$, if there were two independent variables $Y^{(1)},Y^{(2)}$ from the conditional distribution of $Y$ given $X$, then,

\begin{align*}
    Var(1(Y\leq y)|X)=\frac{1}{2} E(1(Y^{(1)}\leq y)-1(Y^{(2)}\leq y)^2)\\
    =\frac{1}{2} E(1(min(Y^{(1)},Y^{(2)})\leq y
    \leq max(Y^{(1)},Y^{(2)}))
\end{align*}

As we have seen in Section 2.4, his final move is to take $Y_{(i)}$ and $Y_{(i+1)}$ as proxies for $Y^{(1)},Y^{(2)}$, and use plug in, leading to $\hat{C}_{n}$.

Note that the same trick can also be applied to estimate $Var\bigg(G_Y(Y)|F_X(X)\bigg)$, needed for $M(X,G_Y(Y))$ which measures directly functional dependence of $Y$ on $X$.

We are led to,

\begin{align}
    \hat{M} \equiv 1-\frac{6}{n} \sum_{i=1}^{n-1} {(r_{i+1}-r_{i})^2}+O(\frac{1}{n}) \notag\\
    =\frac{12}{n} \sum_{i=1}^{n-1} {\frac{r_{i}r_{i+1}}{n^2}}-3+O(\frac{1}{n})\notag\\
    = \hat{T_S}^{(L)}+O(\frac{1}{n})
\end{align}

The corresponding parameter,

\begin{align*}
    M(X,G_Y(Y))=\frac{Var(E(G_Y(Y)|X))}{Var(G_Y(Y))}\\
    =1-\frac{E[Var(G_Y(Y)|X)]}{Var(G_Y(Y))}
\end{align*}

which is just $T_{S}^{(L)}$, the analogue of Spearman's $\rho$.

It satisfies C) (Functional Dependence) but not B).

which can be rewritten as we did $C$ in (2.32) as

\begin{align}
=1-\frac{E_{F^{(2)}}\bigg( F_2(Y^{(1)})-F_2(Y^{(2)})\bigg)^2} {E\bigg( F_2(Y_1)-F_2(Y_2)\bigg)^2}    
\end{align}

This is equal if $F_1$ and $F_2$ are continuous to $T_{S}^{(L)}(F)$ which is consistently estimated by $\hat{T}_{S}^{(L)}$ given in (3.9).

Expression (3.10) suggests a class of statistics which have the same power and functional dependence properties as $T_S^{(L)}(F)$ and $C$.

Given a monotone strictly increasing on $(0,1)$, define the statistic

\begin{align*}
    \hat{T}_{a}^{(L)}=\frac{\frac{1}{n-1}\sum_{i=1}^{n-1}a(\frac{r_i}{n})a(\frac{r_{i+1}}{n})-\bar{a}^2}{\sigma^2(a)}
\end{align*}

where
\begin{align*}
    \bar{a} \equiv \frac{1}{n}\sum_{i=1}^n a(\frac{i}{n})\\
    \sigma^2(a) \equiv \frac{1}{n}\sum_{i=1}^n\bigg(a(\frac{i}{n})-\bar{a}\bigg)^2
\end{align*}

$\hat{T}_{a}^{(L)}$ under mild conditions given below consistently estimates,

\begin{align*}
    T_{a}^{(L)}(F)=1-\frac{E_{F^{(2)}}\bigg( a\big(F_2(Y^{(1)})\big)-a\big(F_2(Y^{(2)})\big)\bigg)^2} {E\bigg( a(F_2(Y_1))-a(F_2(Y_2))\bigg)^2} 
\end{align*}

the analogues of objects each as the normal scores correlation (if $a=\Phi^{-1}$) 

We prove,

\begin{theorem}
Assume $(U,V)\in [0,1]^2$ are such that $U \sim \mc{U}(0,1)$, $V$ given $U=u$ has cdf $F(.|u)$ and $V$ has a continuous cdf,

\begin{align*}
     &(i)\int|dF(v|u)-dF(v|u_0)| \to 0 \textit{ whenever } u \to u_0\\
    &\textit{(i) can be interpreted as continuity of the map $u \to F(.|u)$ in $\mc{P}_0$ variational norm on probabilities.}\\
    &(ii) q: [0,1]^2
 \to R, q \textit{ continuous}.
\end{align*}

Then,
\begin{align}
    \frac{1}{n-1}\sum_{i=1}^{n-1}q\big(\frac{r_i}{n},\frac{r_{i+1}}{n}\big) \stackrel{p}{\lra} \int_0^1 \int_0^1 \int_0^1 q(v_1,v_2)dF(v_1|u)dF(v_2|u)du
\end{align}
\end{theorem}

\textbf{Proof of Theorem 3.1}

Write,
\begin{align*}
    \frac{1}{n-1}\sum_{i=1}^{n-1}q\big(\frac{r_i}{n},\frac{r_{i+1}}{n}\big)=\frac{1}{n-1}\sum_{i=1}^{n-1}q\big(\hat{F}_2(Y_{[i]}),\hat{F}_2(Y_{[i+1]})\big) 
\end{align*}
 
where $Y_{[i]}=Y_j$ iff $X_j=X_{(i)}$,
where $X_{(i)}$ are the ordered $X_i$ as usual. We can write $X_{(i)} \equiv F_1^{-1}(U_{(i)})$, $Y_{(i)} \equiv F_2^{-1}(V_{[i]})$ with $(U_i,V_i)$ defined appropriately by assumption $F_1,F_2$ are continuous.

Note that given $X_{(1)},...,X_{(i)}$, 
the $Y_{[i]}$ are independent with $Y_{[i]}\sim \mathcal{L}(Y|X=X_{(i)})$.

Next use Glivenko-Cantelli on the $Y_1,...,Y_n$ and the uniform continuity of $q$ to conclude that

\begin{align*}
    \frac{1}{n-1}\sum_{i=1}^{n-1}q(r_i,r_{i+1})-\frac{1}{n-1}\sum_{i=1}^{n-1}q\big( V_{[i]},V_{[i+1]}\big) \stackrel{a.s.}{\lra} 0
\end{align*}

But $V_{[1]},...,V_{[n]}$ is a double array of bounded conditionally independent random variables.

Therefore, if $n=2M+2$, by the LLN for double arrays,

\begin{align}
     \frac{1}{n-1}\sum_{k=0}^{n} \bigg[q\big( V_{[2k+1]},V_{[2k+2]}\big)-E\big[ q\big( V_{[2k+1]},V_{[2k+2]}\big)|U_{(2k+1)},U_{(2k+2)}\big]\bigg]\stackrel{p}{\lra} 0
\end{align}

Using $(V_{[2k]},V_{[2k+1]})$ in the same way we used $(V_{[2k+1]},V_{[2k+2]})$ in (3.12) we obtain an analogous result.

Let
\begin{align*}
    Q(u_1,u_2) \equiv E\bigg[q(V_1,V_2)|U_1=u_1,U_2=u_2\bigg]
\end{align*}

If we add the first sum in (3.12) to the corresponding sum with $(V_{[2k]},V_{[2k+1]})$, we see that the theorem is proved if we can show that

\begin{align}
    \frac{1}{n-1}\sum_{i=1}^{n-1}Q(U_{(i)},U_{(i+1)})= \frac{1}{n}\sum_{i=1}^{n} Q(U_{(i)},U_{(i)})+o_p(1)
\end{align}

Then we can just apply the LLN to 

\begin{align*}
    \frac{1}{n}\sum_{i=1}^{n}Q(U_{(i)},U_{(i)})=\frac{1}{n}\sum_{i=1}^{n}Q(U_i,U_i)
\end{align*}

To prove (3.13) we write,

\begin{align*}
    Q(x_1,x_2)=\int q(v_1,v_2)dF(v_1|x_1)dF(v_2|x_2)
\end{align*}

Then,

\begin{align*}
    Q(X_{(i)},X_{(i+1)})-Q(X_{(i)},X_{(i)})=\int \int q(v_1,v_2)dF(v_1|X_{(i)})\big(dF(v_2|X_{(i+1)})-dF(v_2|X_{(i)})\big)
\end{align*}

We can now apply condition (i) using the boundedness of $q$ guaranteed by (ii) to establish (3.13). \qedsymbol

\textbf{Remarks: }

1) Since $a(.)$ in $T_a^{(L)}$ is assumed continuous, consistency of $T_a^{(L)}$ follows from Theorem 3.1, as does consistency of $\hat{C}_n$.

2) Condition (i) can be weakened using the methods of Chatterjee \cite{ChatterjeeNCCA2020}.

\section{The best of both worlds}

We will now argue that decoupling independence and functional dependence is the "right" thing to do. We want to construct measures of dependence $D(F)$ and corresponding empirical versions $\hat{D}_n$, such that, $\hat{D}_n \to D(F)$ in probability for all cdf $F$ of $(X,Y)$ and these quantities have the following properties:

\begin{enumerate}[label=(\roman*)]
    \item $0 \leq D(F) \leq 1$
    \item $D(F)=0$ iff $X$ and $Y$ are independent
    \item $D(F)=1$ iff $Y=h(X)$ a.s. for some $h$
    \item $\hat{D}_n$ based tests have "optimal" local power properties against $H$: $X$ and $Y$ are independent.
\end{enumerate}

Let $T_I(F)$ be any measure of dependence such that $0 \leq T_I(F) \leq 1$, $T_I(F) =0$ iff $X$ and $Y$ are independent. Let $\hat{T_I}$ be an associated statistic which is consistent whatever be $F$, and has, as a test statistic, any desired local power properties. From our discussion this might be power against any local alternative for which $\dot{l}_0(X,Y) \neq a(X)+b(Y)$, $a,b \neq 0$. Or it might be particularly good local power behaviour in some directions as in Bickel, Ritov and, Stoker (\cite{BickelTTGFA2006}). Suppose also that $T_{I}(F)<1$, for all $F$. This parameter may now naturally be coupled to $T_{FD}(F)$ such that $T_{FD}(F)=1$ iff $Y=h(X)$ a.s., $\hat{T_{FD}}$, is consistent for $T_{FD}(F)$, and $H \Rightarrow T_{FD}(F)=0$.

Let,
\begin{align*}
    D(F)&=T_I(F) \textit{ if } 0 \leq T_{FD}(F) \leq \frac{1}{2} \\
    &= T_{FD}(F) \textit{ otherwise.}
\end{align*}

Let $\hat{D}$ be similarly defined in terms of $\hat{T_I}$ and $T_{FD}(\hat{F})$. 

Note that $\frac{1}{2}$ may be replaced by any number $<1$.

Then, clearly,

\begin{align*}
    &(i) P_0[\hat{D}=\hat{T_I}] \to 1 \\
    &(ii) \hat{D} \stackrel{p}{\lra} 0 \Rightarrow H\\
    &(iii) \hat{D} \stackrel{p}{\lra} 1\\
    &\textit{iff } Y=h(X) \textit{ a.s.} (FD)
\end{align*}

(i)-(iii) are argued as follows.

(i) $H\Rightarrow \hat{T}_I \stackrel{p}{\lra} 0$ and $\hat{T}_{FD} \stackrel{p}{\lra} 0 \Rightarrow \hat{D} \stackrel{p}{\lra} 0$.

(ii)$\hat{D} \stackrel{p}{\lra} 0 \Rightarrow \hat{T}_I \stackrel{p}{\lra} 0$ AND $\hat{T}_{FD} \stackrel{p}{\lra} 0$.

The first implication implies $H$.

The second implies $P[\hat{D}=\hat{T}_I] \to 1$ by definition, so that the first case holds as well.

(iii)  a)$FD \Rightarrow \hat{T}_{FD} \stackrel{p}{\lra} 1 \Rightarrow P[\hat{D}=\hat{T}_{FD}] \to 1 \Rightarrow \hat{D} \stackrel{p}{\lra} 1$, 

\quad b) $\hat{D} \stackrel{p}{\lra} 1 \Rightarrow P[\hat{D}=\hat{T}_{FD}] \to 1 \Rightarrow FD$

Both a) and b) hold by consistency and the assumption that $T_{I}(F)<1$.

Thus $D, \hat{D}$ have the properties we want. For instance, we can take $\hat{T}_{I}= \hat{T}_{BKR}$.

and $\hat{T}_{FD}=\hat{C}_n$ or $\hat{T}_{S}^{(L)}$. Or we could add desired local power properties by taking 

\begin{align*}
\hat{T}_{I}=max\{ \hat{T}_{BKR},\hat{T}_{BKR}^{(L)} \}.    
\end{align*}

We use Lemma A2, establishing $T_{BKR}^{(F)}<1$. The same property applies to $T^{(L)}_{BKR}$ by the same argument. Clearly, $H \Rightarrow \hat{T}_S^{(L)} \stackrel{p}{\lra} 0$ and FD holds iff $\hat{T}_S^{(L)} \stackrel{p}{\lra} 1$.

$\bold{Discussion:}$ As Lin and Han (2021) show, it is possible by going to parameters which depend on $F^{(m)},m \to \infty$ to obtain extensions of $\hat{C_n}$ which have good minimax properties beyond the $n^{-\frac{1}{4}}$ detection barrier, as they define it.

This could be done more generally with parameters like $D(F)$ but we do not pursue this or generalizations to other spaces than $R$, as done by Deb, Ghosal and Sen (2020).

\begin{appendices}

\section{Lemma A.1}

\begin{lemma}
 If $\dot{l}_0 \neq 0$ satisfies the hypothesis of Theorem 2.1, then $C(.)$ given in (2.27) is not 0.
\end{lemma}

$\bold{Proof:}$

By definition,

\begin{align*}
    T_n(.) \stackrel{d}{\lra} T(.)+tC(.)
\end{align*}

where,

\begin{align}
    C(\tend{u})=\int \Psi(\tend{u},\tend{w}) \dot{l}_0(\tend{w}) d \tend{w}
\end{align}

and

\begin{align}
   \Psi(\tend{u},\tend{w})=1(\tend{w} \leq \tend{u})- w_1 1(w_2 \leq u_2) -w_2 1(w_1 \leq u_1)+w_1w_2 
\end{align}

Here, $\tend{u}=(u_1,u_2)$, $\tend{w}=(w_1,w_2)$, $\tend{w} \leq \tend{u}$ iff $w_j \leq u_j, j=1,2$.

Rewriting,

\begin{align}
    C(\tend{u})=E\Psi(\tend{u},U,V) \dot{l}_0(U,V)
\end{align}

for $U,V$ iid $\mathcal{U}(0,1)$. By construction, for all $\tend{u}$,

\begin{align}
    E(\Psi(\tend{u},U,V)|U)=E(\Psi(\tend{u},U,V)|V)=0
\end{align}

Therefore, if $C(\tend{u})=0$ for all $\tend{u}$, $Ef(U,V)\dot{l}_0(U,V) =0$ for all $f \in \Lambda_1^{\perp}$.

But,  $ \dot{l}_0 \in \Lambda_1^{\perp}$.

The lemma follows. \qedsymbol

\begin{lemma}
 $T_{BKR}(F)<1$ for all $F$.
\end{lemma}

$\bold{Proof:}$ Suppose $F_1,F_2$ are continuous. Then $T_{BKR}(F)$ is given by (2.2).

By (2.2) we an take $F_j(u)=u$ for $j=1,2$. 

Since $|F(u,v)-uv| \leq 1$ for all $u,v$, we must have,

\begin{align*}
    |F_m(u,v)-uv| \to 1 
\end{align*}

in measure on $[0,1]^2$. Since the sequence $\{F_m\}$ is tight, there exists $\{m_k\}$, $F_0$ such that $F_{m_k} \stackrel{d}{\lra} F_0$ and

\begin{align*}
    |F_0(u,v)-uv|=1 \textit{ a.e. } (u,v)
\end{align*}

equivalent to,

\begin{align*}
    F_0(u,v)=1+uv \textit{ or } uv-1.
\end{align*}

This is impossible in the interior of $[0,1]^2$. The result follows for $F_1, F_2$ continuous.

In general, write,

\begin{align*}
    T_{BKR}(F)=\int_0^1\int_0^1 \bigg( F\big(F_1^{-1}(u),F_2^{-1}(v)\big) -F_1\big(F_1^{-1}(u)\big) F_2\big(F_2^{-1}(v)\big)  \bigg)^2 dudv
\end{align*}

The argument given before holds if $T_{BKR}=1$,

\begin{align}
    F_1\big(F_1^{-1}(u)\big) F_2\big(F_2^{-1}(v)\big)=0 \textit{ and } F\big(F_1^{-1}(u),F_2^{-1}(v)\big)=1
\end{align}

or conversely for almost all $u,v$.

But if $(X,Y)$ takes on only 1 value,

\begin{align*}
    F_1F_2=F \textit{ a.e.}
\end{align*}

While if $(X,Y)$ takes on more than 1 value but is not continuous, (A.5) is impossible since $F(.,.)$ takes on more than 1 value. In either case,m the lemma holds.
\qedsymbol

\section{Proof of Lemma 2.2}

We begin with an identity based on Fubini's theorem. For any finite signed measure $\mu$ on $[0,1]^2$ if $\tilde{A}(s)=A(1)-A(s)$, $\tilde{B}$ is similarly defined, and $\mu(s,t)=\int_0^s\int_0^td\mu(u,v)$, then

\begin{align}
    \int_0^1\int_0^1 \bar{\tilde{A}}(s) \bar{\tilde{B}}(t) d\mu(s,t)=\int_0^1\int_0^1 a(s)b(t) \mu(s,t) dsdt
\end{align}

Note that, $\bar{\tilde{A}}(s)=-\bar{A}(s)$.

Then, if $C(u,v)=uv$,

\begin{align}
   \frac{1}{\sqrt{n}} \sum_{k=1}^n \bar{A}(\frac{R_k}{n})\bar{B}(\frac{S_k}{n})=\sqrt{n}\int_0^1\int_0^1 \bar{\tilde{A}}(u) \bar{\tilde{B}}(v) d(\hat{C}-C)(u,v)
\end{align}

since $\int_0^1\bar{A}(u)du=\int_0^1\bar{B}(v)dv=0$.

Apply (A.5) with $\mu=\hat{C}-C$ to get,

\begin{align*}
\int_0^1\int_0^1 \bar{\tilde{A}}(u) \bar{\tilde{B}}(v) d(\hat{C}-C)(u,v)=\int_0^1\int_0^1 a(s)b(t)\sqrt{n}(\hat{C}-C)(s,t)dsdt
\end{align*}

Apply Proposition 2.1 to get,

\begin{align}
    \int_0^1\int_0^1 a(s)b(t)\sqrt{n}(\hat{C}-C)(s,t)dsdt=\int_0^1\int_0^1a(s)b(t) (E_n(s,t)-sE_{n2}(t)-tE_{n1}(s))dsdt+o_p(1)
\end{align}

since $sup\{|E_n(s,t)|\}=O_p(1)$.

Apply (B.1) again to $E_n(s,t)-sE_{n2}(t)-tE_{n1}(s)$, which is of bounded variation for each n to finally obtain,

\begin{align*}
    &\int_0^1\int_0^1a(s)b(t) (E_n(s,t)-sE_{n2}(t)-tE_{n1}(s))dsdt=\\
    &=\int_0^1\int_0^1 \bar{\tilde{A}}(s) \bar{\tilde{B}}(t)dE_n(s,t)
\end{align*}

since $\int_0^1ta(t)dt=\int_0^1\bar{A}(t)dt$. The lemma follows.
\qedsymbol

\section{Proof of Lemma 2.3}

Let $B_n= \sigma$ field generated by $U_1,...,U_n$.

We can then write $A_n=A_{n1}+A_{n2}$, where

\begin{align*}
    A_{n1}\equiv E(A_n|B_n)&=\frac{1}{\sqrt{n}}\sum_{k=1}^n\int_0^1 a(U_k,v)dv\\
    &=\frac{1}{\sqrt{n}}\sum_{k=1}^n a_1(U_k)
\end{align*}

\begin{align*}
    A_{n2}\equiv A_n-A_{n1}
\end{align*}

Given $B_n$, $A_{n2}$ is a sum of n independent, non identically distributed variables with mean 0 and 
\begin{align*}
    Var(A_{n2}|B_n)=\frac{1}{n}\sum_{k=1}^n Var(a(U_k,V_k)|V_k)
\end{align*}

By the Lindeberg Feller theorem and the LOL, 

\begin{align*}
    Var(A_{n2}|B_n) \stackrel{a.s.}{\lra} E\{Var[a(U_1,V_1)|U_1]\} \equiv \sigma_2^2
\end{align*}

and

\begin{align*}
    A_{n2} \stackrel{d}{\lra} Z_0^{(2)} \sim N(0,\sigma_2^2)
\end{align*}

$A_{n1}$ is uncorrelated with $A_{n2}$ and by the CLT,

\begin{align*}
    A_{n1} \stackrel{d}{\lra} Z_0^{(1)} \sim N(0,\sigma_1^2), \sigma_1^2\equiv E a_1^2(U_1)
\end{align*}

Therefore,

\begin{align*}
    (A_{n1},A_{n2}) \stackrel{d}{\lra}(Z_0^{(1)},Z_0^{(2)})
\end{align*}

where $Z_0^{(1)},Z_0^{(2)}$ are independent $N(0,\sigma_1^2)$, $N(0,\sigma_2^2)$.

Consider $h^{(1)}(v_1,v_2)$, and write

\begin{align*}
    H_n^{(1)} = \frac{1}{\sqrt{n}} \sum_{k=1}^{n-1} h^{(1)}(W_k,W_{k+1})
\end{align*}

where $W_k=V_{D_k}$ as before are iid $\mathcal{U}(0,1)$ given $B_n$. By Orey's (1958) \cite{OreyCLTMD1958} theorem for 2 dependent random vairable, conditionally,

\begin{align*}
    H_n^{(1)} \stackrel{d}{\lra} Z_1 \sim N(0, Var(Z_1))
\end{align*}
with $Var(Z_1)$ as specified in (2.34).

The same applies to $H_n^{(2)}$ and the joint distribution is as specified. We can put $(A_{n2},H_n^{(1)},H_n^{(2)})$ together to obtain joint convergence given $B_n$. Since the limiting covariance matrix is independent of $B_n$, we can add $A_{n1}$ as asymptotically independent of $(A_{n2},H_n^{(1)},H_n^{(2)})$. Since the covariances of the limits are in all cases the limit of the unconditional covariances the lemma follows on direct computation of the latter.

\qedsymbol

\section{Proof of Proposition 2.4}

Consider for $(X_i,Y_i)$ iid $h_0$, and apply Taylor's theorem to

\begin{align}
    &\sum_{i=1}^n \big[ l_{\theta}(X_i,Y_i)-l_{0}(X_i,Y_i)\big] \notag\\
    &=\theta \sum_{i=1}^n \big(a(X_i)+b(Y_i)\big)+\frac{\theta^2}{2} \sum_{i=1}^n \frac{\ptl^2 l_{\theta}}{\ptl \theta^2} \bigg|_{0}(X_i,Y_i) +\frac{\theta^3}{6} R_n(\theta)
\end{align}

By assumption, for $|\theta|\leq \frac{M}{\sqrt{n}}$,

\begin{align}
    |R_n(\theta)|\leq \frac{M^3}{6n^{\frac{1}{2}}} K=o(1)
\end{align}

On the other hand, since $Ea(X)=Eb(Y)=0$, $a,b$ bounded, under $H$,

\begin{align}
    \theta_n \sum_{i=1}^n \big(a(X_i)+b(Y_i)\big) \stackrel{d}{\lra} Z
\end{align}

if $\sqrt{n}\theta_n \to t, \sigma^2 \equiv E_0\big(a(X)+b(Y)\big)^2, Z \sim N(0,t^2\sigma^2)$.

By a standard argument, under our conditions,

\begin{align*}
    E_0\big(\frac{\ptl^2}{\ptl \theta^2}l_\theta(X_i,Y_i)\big|_{0}\big)=-\sigma^2
\end{align*}

Thereforem, by the LLN, if $\sqrt{n}\theta_n \to t$,

\begin{align}
    \frac{\theta_n^2}{2} \sum_{i=1}^n \frac{\ptl^2}{\ptl \theta^2}l_\theta(X_i,Y_i)\big|_{0} \stackrel{a.s}{\lra} -\frac{t^2}{2} \sigma^2
\end{align}

The same argument applies using only the boundedness of $a,b$ to $\tilde{h}(X_i,Y_i)$. Write,

\begin{align}
    \Delta_n^{(1)}(\theta) = E_0 \big| \prod_{i=1}^n \frac{h_\theta(X_i,Y_i)}{h_0(X_i,Y_i)} -exp \{\theta \sum_{i=1}^n (a(X_i)+b(Y_i))-\frac{\theta^2}{2} \sigma^2\}  \big|
\end{align}

and define $\Delta_n^{(2)}(\theta)$ similarly for the $\tilde{h}(X_i,Y_i)$.

Then,

\begin{align*}
    \Delta_n(\theta) \leq \Delta_n^{(1)}(\theta)+\Delta_n^{(2)}(\theta)
\end{align*}

We can write, letting $A_n(.)$ and $B_n(.)$ be the logs of the expressions in (D.5),

\begin{align*}
    \Delta_{n}^{(1)}(\theta) \equiv E_0|e^{A_n(\theta)}-e^{B_n(\theta)}|
\end{align*}

By (D.2),

\begin{align}
|A_n(\theta)-B_n(\theta)|=o_p(1)    
\end{align}
 
Moreover, by (D.3) and (D.4),

\begin{align*}
    e^{B_n(\theta)} \stackrel{d}{\lra} e^{A_t}
\end{align*}

where $A_t=Z-\frac{t^2\sigma^2}{2}$ and $Ee^{A_t}=1$.

Since $E_0 e^{A_n(\theta)}=1$ for all $n$, we can conclude by uniform integrability of the $A_n(\theta)$ that,

\begin{align*}
    E_0 e^{A_n(\theta)} |1-e^{B_n(\theta)-A_n(\theta)}| \to 0
\end{align*}

The same argument applies to $\Delta_n^{(2)}(\theta) $ and we can conclude that if $\sqrt{n}\theta_n \to t$,

\begin{align*}
    \Delta_n(\theta_n) \to 0 
\end{align*}

Finally, we note that if $sup\{\Delta_n(\theta): |\theta|\leq \frac{M}{\sqrt{n}}\}$ does not tend to 0, there must exist $\theta_n$ with $|\sqrt{n}\theta_n| \to |t| \leq M$ for which $\Delta_n(\theta_n)$ does not tend to 0, a contradiction. 

\qedsymbol

\end{appendices}

\section*{Acknowledgements}
This work is based in part on joint work with Sky Cao \cite{CaoCTEPO2020}. His contribution to the current paper amounts to coauthorship. However, any errors are the responsibility of the author.

We thank Sourav Chatterjee for facilitating this collaboration, as well as for helpful conversations. We thank Holger Dette for pointing out an important reference. We thank Hongjian Shi and Mathias Drton for comments on the local power calculations. In particular, we thank Fang Han also for valuable comments on various parts of the paper. 

\bibliographystyle{abbrvnat}
\bibliography{peter-old-3}

\end{document}